\documentclass{svmult}

\renewcommand{\email}[1]{\emailname: #1} 

\usepackage{mathptmx}       
\usepackage{helvet}         
\usepackage{courier}        

\usepackage{makeidx}         
\usepackage{graphicx}        
\usepackage[bottom]{footmisc}

\usepackage{latexsym}
\usepackage{amsmath}
\usepackage{amsfonts}
\usepackage{amssymb}
\usepackage{bm}

\usepackage{url}
\usepackage[misc,geometry]{ifsym}

\spdefaulttheorem{assumption}{Assumption}{\upshape \bfseries}{\itshape}
\spdefaulttheorem{algo}{Algorithm}{\upshape \bfseries}{\itshape}

\newcommand{\bsd}{{\boldsymbol{d}}}
\newcommand{\bse}{{\boldsymbol{e}}}

\newcommand{\bsr}{{\boldsymbol{r}}}

\newcommand{\bsx}{{\boldsymbol{x}}}
\newcommand{\bsy}{{\boldsymbol{y}}}
\newcommand{\bsz}{{\boldsymbol{z}}}
\newcommand{\bsA}{{\boldsymbol{A}}}

\newcommand{\bsD}{{\boldsymbol{D}}}

\newcommand{\bsI}{{\boldsymbol{I}}}

\newcommand{\bsX}{{\boldsymbol{X}}}
\newcommand{\bsY}{{\boldsymbol{Y}}}
\newcommand{\bsZ}{{\boldsymbol{Z}}}

\newcommand{\bszero}{{\boldsymbol{0}}} 

\newcommand{\bsphi}{{\boldsymbol{\phi}}}

\DeclareSymbolFont{bbold}{U}{bbold}{m}{n}
\DeclareSymbolFontAlphabet{\mathbbold}{bbold}

\newcommand{\ceil}[1]{\left\lceil #1 \right\rceil}    

\pdfoutput=1

\begin{document}

\title*{Efficient Spherical Designs with Good Geometric Properties}
\titlerunning{Efficient Spherical Designs}

\author{Robert S. Womersley}

\institute{
 Robert S. Womersley (\Letter)
 \at School of Mathematics and Statistics, University of New South Wales, Sydney NSW 2052 Australia \\
 \email{r.womersley@unsw.edu.au}
}

\maketitle

\index{Womersley, Robert S.}

\paragraph{Dedicated to Ian H.~Sloan on the occasion of his $80$th birthday
in acknowledgement of his many fruitful ideas and generosity.}

\abstract{
Spherical $t$-designs on $\mathbb{S}^{d}\subset\mathbb{R}^{d+1}$ provide $N$ nodes for
an equal weight numerical integration rule which is exact for all spherical polynomials
of degree at most $t$. This paper considers the generation of efficient,
where $N$ is comparable to $(1+t)^d/d$, spherical $t$-designs with
good geometric properties as measured by their mesh ratio,
the ratio of the covering radius to the packing radius.
Results for $\mathbb{S}^{2}$ include computed spherical $t$-designs for $t = 1,...,180$
and symmetric (antipodal) $t$-designs for degrees up to $325$, all with low mesh ratios.
These point sets provide excellent points for numerical integration on the sphere.
The methods can also be used to computationally explore spherical $t$-designs for $d = 3$
and higher.
}

\section{Introduction}\label{RSWsec:Intro}

Consider the $d$-dimensional unit sphere
\[
\mathbb{S}^{d} = \left\{ \bsx\in \mathbb{R}^{d+1}: | \bsx | = 1 \right\}
\]
where the standard Euclidean inner product is
$\bsx \cdot \bsy = \sum_{i=1}^{d+1} x_i y_i$ 
and $| \bsx |^2 = \bsx \cdot \bsx$.

A numerical integration (quadrature) rule for $\mathbb{S}^{d}$ is a set of
$N$ points $\bsx_j \in \mathbb{S}^{d}, j = 1,\ldots,N$ and associated weights $w_j > 0, j = 1,\ldots,N$
such that
\begin{equation}
  Q_N(f) := \sum_{j=1}^N w_j f(\bsx_j) \approx
  I(f) := \int_{\mathbb{S}^{d}} f(\bsx) d \sigma_d(\bsx) .
  \label{RSWeq:QN}
\end{equation}
Here $\sigma_d(\bsx)$ is the normalised Lebesgue measure on $\mathbb{S}^{d}$
with surface area
\[
  \omega_d := 
  \frac{2 \pi^{(d+1)/2}}{\Gamma((d+1)/2)},
\]
where $\Gamma(\cdot)$ is the gamma function.

Let $\mathbb{P}_t(\mathbb{S}^{d})$ denote the set of all spherical polynomials on $\mathbb{S}^{d}$ of degree at most $t$.
A \emph{spherical $t$-design} is a set of $N$ points $X_N=\{\bsx_1,\ldots,\bsx_N\}$ on $\mathbb{S}^{d}$
such that equal weight quadrature using these nodes is exact for all spherical polynomials
of degree at most $t$, that is
\begin{equation}
  \frac{1}{N} \sum_{j=1}^N p(\bsx_j) =
  \int_{\mathbb{S}^{d}} p(\bsx) \mbox{d}\sigma_d(\bsx),
  \quad \forall p \in \mathbb{P}_t(\mathbb{S}^{d}).
\label{RSWeq:sd}
\end{equation}

Spherical $t$-designs were introduced by Delsarte, Goethals and Seidel~\cite{DelGS77}
who provided several characterizations and established lower bounds on the number of points
$N$ required for a spherical $t$-design.
Seymour and Zaslavsky\cite{SeyZla84} showed that spherical $t$-designs exist on $\mathbb{S}^{d}$
for all $N$ sufficiently large.
Bondarenko, Radchenko and Viazovska~\cite{BonRV13} established that there exists a $C_d$
such that spherical $t$-designs on $\mathbb{S}^{d}$ exist for all $N \geq C_d \; t^d$,
which is the optimal order.
The papers \cite{ConSlo99,BanBan09,CohKum07} provide a sample of many
on spherical designs and algebraic combinatorics on spheres.

An alternative approach, not investigated in this paper,
is to relax the condition $w_j = 1/N$ that the quadrature weights are equal
so that $|w_j / (1/N) - 1| \leq \epsilon$ for $j = 1,\ldots,N$ and $0 \leq \epsilon < 1$,
but keeping the condition that the quadrature rule is exact for polynomials of degree $t$
(see \cite{SloWom04,ZhoChe17} for example).

The aim of this paper is not to find spherical $t$-designs with the
minimal number of points, nor to provide proofs that a particular configuration
is a spherical $t$-design. Rather the aim is to find sequences of point sets
which are at least computationally spherical $t$-designs,
have a low number of points and are geometrically well-distributed on the sphere.
Such point sets provide excellent nodes for numerical integration on the
sphere, as well as hyperinterpolation~\cite{Sloan95,HesSlo07,SloWom12}
and fully discrete needlet approximation~\cite{WanLSW16}.
These methods have a requirement that the quadrature rules are
exact for certain degree polynomials.
More generally, \cite{HesSW10} provides a summary of numerical integration on $\mathbb{S}^{2}$
with geomathematical applications in mind.

\subsection{Spherical Harmonics and Jacobi Polynomials}
A \emph{spherical harmonic} of degree $\ell$ on $\mathbb{S}^{d}$ is the restriction to $\mathbb{S}^{d}$
of a homogeneous and harmonic polynomial of total degree $\ell$ defined on $\mathbb{R}^{d+1}$.
Let $\mathbb{H}_{\ell}$ denote the set of all spherical harmonics of exact degree $\ell$ on $\mathbb{S}^{d}$.
The dimension of the linear space $\mathbb{H}_{\ell}$ is
\begin{equation}\label{eq:dim.sph.harmon}
    Z(d,\ell):=(2\ell+d-1)\frac{\Gamma(\ell+d-1)}{\Gamma(d)\Gamma(\ell+1)}\asymp (\ell+1)^{d-1},
\end{equation}
where $a_{\ell}\asymp b_{\ell}$ means $c\:b_{\ell}\le a_{\ell}\le c' \:b_{\ell}$ for some
positive constants $c$, $c'$, and the asymptotic estimate uses \cite[Eq.~5.11.12]{NIST:DLMF}.

Each pair $\mathbb{H}_{\ell}$, $\mathbb{H}_{\ell'}$ for $\ell\neq \ell'\ge0$ is
$\mathbb{L}_{2}$-orthogonal, $\mathbb{P}_{L}(\mathbb{S}^{d})=\bigoplus_{\ell=0}^{L} \mathbb{H}_{\ell}$
and the infinite direct sum $\bigoplus_{\ell=0}^{\infty} \mathbb{H}_{\ell}$
is dense in $\mathbb{L}_p(\mathbb{S}^d)$, $p \geq 2$, see e.g. \cite[Ch.1]{WaLi2006}.
The linear span of $\mathbb{H}_{\ell}$, $\ell=0,1,\dots,L$, forms the space
$\mathbb{P}_{L}(\mathbb{S}^{d})$ of spherical polynomials of degree at most $L$.
The dimension of $\mathbb{P}_{L}(\mathbb{S}^{d})$  is
\begin{equation}
  D(d, L) := \mbox{dim}\; \mathbb{P}_{L}(\mathbb{S}^{d}) = \sum_{\ell=0}^L Z(d, \ell) = Z(d+1,L).
  \label{RSWeq:dimP}
\end{equation}

Let
$P^{\left(\alpha,\beta\right)}_{\ell}(z)$, $-1\le z\le1$,
be the Jacobi polynomial of degree $\ell$ for $\alpha,\beta>-1$.
The Jacobi polynomials form an orthogonal polynomial system with respect to the Jacobi weight
$w_{\alpha,\beta}(z) := (1-z)^{\alpha}(1+z)^{\beta}$, $-1\le z \le 1$.
We denote the normalised Legendre (or ultraspherical/Gegenbauer) polynomials by
\begin{equation*}
  P^{(d+1)}_{\ell}(z) := \frac{P^{\left(\frac{d-2}{2},\frac{d-2}{2}\right)}_{\ell}(z)}{P^{\left(\frac{d-2}{2},\frac{d-2}{2}\right)}_{\ell}(1)},
\end{equation*}
where, from \cite[(4.1.1)]{Szego75},
\begin{equation}
  P_\ell^{(\alpha,\beta)}(1) = \frac{\Gamma(\ell+\alpha+1)}{\Gamma(\ell+1)\Gamma(\alpha+1)},
  \label{RSWeq:Pab1}
\end{equation}
and \cite[Theorem~7.32.2, p.~168]{Szego75},
\begin{equation}
  \bigl|P^{(d+1)}_{\ell}(z)\bigr|\le 1, \qquad -1\leq z \leq 1.
  \label{RSWeq:LegBnd}
\end{equation}
The derivative of the Jacobi polynomial satisfies \cite{Szego75}
\begin{equation}
  \frac{\mbox{d}\; P_{\ell}^{(\alpha,\beta)}(z)}{\mbox{d}\; z} =
  \frac{\ell+\alpha+\beta+1}{2}\; P_{\ell-1}^{(\alpha+1, \beta+1)}(z),
\end{equation}
so
\begin{equation}
  \frac{\mbox{d}\; P_{\ell}^{(d+1)}(z)}{\mbox{d}\; z} =
  \frac{(\ell+d-1)(\ell + d/2)}{d} P_{\ell-1}^{(d+3)}(z).
\end{equation}
Also if $\ell$ is odd then the polynomials $P_\ell^{(d+1)}$ are odd
and if $\ell$ is even the polynomials $P_\ell^{(d+1)}$ are even.

A \emph{zonal function} $K:\mathbb{S}^{d}\times\mathbb{S}^{d}\rightarrow \mathbb{R}$
depends only on the inner product of the arguments,
i.e. $K(\bsx,\bsy)= \mathfrak{K}(\bsx\cdot\bsy)$,\: $\bsx,\bsy\in \mathbb{S}^{d}$,
for some function $\mathfrak{K}:[-1,1]\to \mathbb{R}$.
Frequent use is made of the zonal function $P^{(d+1)}_{\ell}(\bsx\cdot\bsy)$.

Let $\{Y_{\ell,k}: k=1,\dots,Z(d,\ell), \; \ell = 0,\ldots,L\}$ be an \emph{orthonormal basis} for $\mathbb{P}_L(\mathbb{S}^{d})$.
The normalised Legendre polynomial $P^{(d+1)}_{\ell}(\bsx\cdot\bsy)$
satisfies the \emph{addition theorem}
(see \cite{Szego75,WaLi2006,AtkHan2012} for example)
\begin{equation}\label{RSWeq:addthm}
  \sum_{k=1}^{Z(d,\ell)} Y_{\ell,k}(\bsx) Y_{\ell,k}(\bsy) =
  Z(d,\ell) P^{(d+1)}_{\ell}(\bsx \cdot \bsy).
\end{equation}

\subsection{Number of Points}\label{RSWsec:Npts}

Delsarte, Goethals and Seidel~\cite{DelGS77} showed that an $N$ point $t$-design on $\mathbb{S}^{d}$ has
$N\ge N^*(d,t)$ where
\begin{equation}
  N^*(d,t) := \left\{ \begin{array}{ll}
       2\binom{d+k}{d} &  \text{ if } t = 2k+1,\\[2ex]
       \binom{d+k}{d} + \binom{d+k-1}{d} & \text { if } t = 2k.
       \end{array}\right.
       \label{RSWeq:DGSbnd}
\end{equation}
On $\mathbb{S}^{2}$
  \begin{equation}
 N^*(2,t) := \left\{ \begin{array}{ll}
     \frac{(t+1)(t+3)}{4} &  \text{ if $t$ odd},\\[2ex]
     \frac{(t+2)^2}{4} & \text { if $t$ even}.
     \end{array}\right.
     \label{RSWeq:DGSbndS2}
\end{equation}
Bannai and Damerell~\cite{BanDam79,BanDam80} showed that
\emph{tight spherical $t$-designs} which achieve the lower bounds (\ref{RSWeq:DGSbnd})
cannot exist except for a few special cases (for example except for $t = 1, 2, 3, 5$ on $\mathbb{S}^{2}$)  .

Yudin~\cite{Yudin97} improved
(except for some small values of $d,t$, see Table~\ref{RSWtab:SF1}),
the lower bounds (\ref{RSWeq:DGSbnd}),
by an exponential factor $(4/e)^{d+1}$ as $t\to\infty$,
so $N \geq N^+(d,t) $ where
\begin{equation}
   N^+(d,t) :=
   2 \frac{\int_0^1 (1-z^2)^{(d-2)/2}\;\mbox{d} z}{\int_\gamma^1 (1-z^2)^{(d-2)/2}\;\mbox{d} z} =
   \frac{\sqrt{\pi}\Gamma(d/2)/\Gamma((d+1)/2)}{\int_\gamma^1 (1-z^2)^{(d-2)/2}\;\mbox{d} z},
  \label{RSWeq:Yudinbnd}
\end{equation}
and $\gamma$ is the largest zero of the derivative $\frac{d P_{t}^{(d+1)}(z)}{d z}$
and hence the largest zero of $P_{t-1}^{(\alpha+1,\alpha+1)}(z)$ where $\alpha = (d-2)/2$.
Bounds~\cite{Szego75,AreDGR04} on the largest zero of $P_n^{(\alpha,\alpha)}(z)$ are
\begin{equation}
  \cos\left(\frac{j_0(\nu)}{n+\alpha+1/2}\right) \leq \gamma \leq
  \sqrt{\frac{(n-1)(n+2\alpha-1)}{(n+\alpha-3/2)/(n+\alpha-1/2)}} \;
  \cos\left(\frac{\pi}{n+1}\right),
\end{equation}
where $j_0(\nu)$ is the first positive zero of the Bessel function $J_\nu(x)$.

Numerically there is strong evidence that spherical $t$-designs with
$N = D(2,t) = (t+1)^2$ points exist,
\cite{CheWom06} and \cite{CheFL11} used interval methods to \emph{prove} existence of
spherical $t$-designs with $N = (t+1)^2$ for all values of $t$ up to $100$,
but there is no proof yet that spherical $t$-designs with
$N \leq D(2,t)$ points exist for all degrees $t$.
Hardin and Sloane~\cite{HarSlo96},\cite{HarSloWWW} provide
tables of designs with modest numbers of points,
exploiting icosahedral symmetry.
They conjecture that for $d = 2$ spherical $t$-designs exist
with $N = t^2/2 + o(t^2)$ for all $t$.
The numerical experiments reported here and available from \cite{Womersley_ssd_URL}
strongly support this conjecture.

McLaren~\cite{McLaren63} defined efficiency $E$ for a quadrature rule as
the ratio of the number of independent functions for which the rule is exact
to the number of arbitrary constants in the rule.
For a spherical $t$-design with $N$ points on $\mathbb{S}^{d}$ (and equal weights) 
\begin{equation}
  E = \frac{\mbox{dim} \; \mathbb{P}_t(\mathbb{S}^{d})}{d N} =
  \frac{D(d,t)}{d N}.
  \label{RSWeq:eff}
\end{equation}
In these terms the aim is to find spherical $t$-designs with $E \geq 1$.
McLaren~\cite{McLaren63} exploits symmetry (in particular octahedral and icosahedral)
to seek rules with optimal efficiency.
The aim here is not to maximise efficiency by finding the minimal number of points for a $t$-design on $\mathbb{S}^{d}$,
but rather a sequence of \emph{efficient} $t$-designs with
 $N \asymp \frac{D(d,t)}{d} \asymp \frac{(1+t)^d}{d}$.
Such efficient $t$-designs provide a practical tool for numerical integration and approximation.

\subsection{Geometric Quality}\label{RSWsec:Geom}

The Geodesic distance between two points $\bsx, \bsy \in \mathbb{S}^{d}$ is
\[
  \mbox{dist}(\bsx, \bsy) = \cos^{-1}(\bsx \cdot \bsy),
\]
while the Euclidean distance is
\[
  | \bsx - \bsy| = \sqrt{2(1-\bsx\cdot\bsy)} = 2 \sin(\mbox{dist}(\bsx,\bsy)/2).
\]
The \emph{spherical cap} with centre $\bsz\in\mathbb{S}^{d}$ and radius $\eta \in [0, \pi]$ is
\[
  \mathcal{C}\left(\bsz;\eta\right) = 
  \left\{\bsx\in \mathbb{S}^{d}: \mbox{dist}(\bsx,\bsz) \leq \eta \right\} .
\]
The \emph{separation distance}
\[
  \displaystyle \delta(X_N) = \min_{i\neq j} \; \mbox{dist}(\bsx_i,\bsx_j)
\]
is twice the packing radius for spherical caps of the same radius and centers in $X_N$.
The best packing problem (or Tammes problem) has a long history~\cite{ConSlo99},
starting with \cite{Tammes30,Rankin55}.
A sequence of point sets $\{X_N\}$ with $N \to\infty$ has the optimal order separation
if there exists a constant $c_d^{\textrm{pck}}$ independent of $N$ such that
\[
  \delta(X_N) \geq c_d^{\textrm{pck}} \; N^{-1/d}.
\]

The separation, and all the zonal functions considered in subsequent sections,
are determined by the set of inner products
\begin{equation}
  \mathcal{A}(X_N) := \left\{ \bsx_i \cdot \bsx_j, i = 1,\ldots,N, j = i+1,\ldots,N \right\}
  \label{RSWeq:IPset}
\end{equation}
which has been widely used in the study of spherical codes, see \cite{ConSlo99} for example.
Then
\[
  \max_{z\in\mathcal{A}(X_N)} z = \cos(\delta(X_N)).
\]
Point sets are only considered different if the corresponding sets (\ref{RSWeq:IPset}) differ,
as they are invariant under an orthogonal transformation (rotation) of the point set and
permutation (relabelling) of the points.

The \emph{mesh norm} (or fill radius)
\[
  h(X_N) = \max_{\bsx\in\mathbb{S}^{d}} \; \min_{j=1,\ldots,N} \; \mbox{dist}(\bsx,\bsx_j)
\]
gives the \emph{covering radius} for covering the sphere with spherical caps of the
same radius and centers in $X_N$.
A sequence of point sets $\{X_N\}$ with $N\to\infty$ has the optimal order covering
if there exists a constant $c_d^{\textrm{cov}}$ independent of $N$ such that
\[
  h(X_N) \leq c_d^{\textrm{cov}} \; N^{-1/d}.
\]

The \emph{mesh ratio} is
\[
  \rho(X_N) = \frac{2 h_{X_N}}{\delta_{X_N}} \geq 1.
\]
A common assumption in numerical methods is that the mesh ratio is uniformly bounded,
that is the point sets are \emph{quasi-uniform}.
Minimal Riesz $s$-energy and best packing points can also produce quasi-uniform point sets~\cite{DamMay05,HarSW12,BonHS14}.

Yudin~\cite{Yudin95} showed that a spherical $t$-design with $N$ points
has a covering radius of the optimal order $1/t$.
Reimer extended this to quadrature rules exact for polynomials of degree $t$
with positive weights.
Thus a spherical $t$-design with $N = O(t^d)$ points provides an optimal order covering.

The union of two spherical $t$-designs with $N$ points is a spherical $t$-design with $2N$ points.
A spherical design with arbitrarily small separation can be obtained as
one $N$ point set is rotated relative to the other.
Thus an assumption on the separation of the points of a spherical design
is used to derive results, see \cite{HesLeo08} for example.
This simple argument is not possible if $N$ is less than twice a lower bound
(\ref{RSWeq:DGSbnd}) or (\ref{RSWeq:Yudinbnd}) on the number of points in a spherical $t$-design.

Bondarenko, Radchenko and Viazovska~\cite{BonRV15} have shown that
on $\mathbb{S}^{d}$ well-separated spherical $t$-designs exist for $N \geq c_d^\prime \; t^d$.
This combined with Yudin's result on the covering radius of spherical designs
mean that there exist spherical $t$-designs with $N = O(t^d)$ points
and uniformly bounded mesh ratio.

There are many other ``geometric'' properties that could be used,
for example the spherical cap discrepancy, see \cite{GraTic93} for example,
(using normalised surface measure so $|\mathbb{S}^{d}| = 1$)
\[
  \sup_{\bsx\in\mathbb{S}^{d}, \eta\in [0,\pi]}
  \left| |\mathcal{C}(\bsx,\eta)| -
  \frac{|X_N\cap \mathcal{C}(\bsx,\eta)| }{N}\right|,
\]
or a Riesz $s$-energy, see \cite{BoyDHSS15} for example,
\[
  E_s(X_N) = \sum_{1\leq i < j \leq N} \frac{1}{|\bsx_i - \bsx_j|^s}.
\]

In distinguishing between spherical $t$-designs with  the
same number $N$ of points we prefer those with lower mesh ratio.
Note that some authors, see \cite{HarSW12,BonHS14} for example,
define the mesh ratio as
$\tilde\rho(X_N) = h(X_N)/\delta(X_N) \geq 1/2$.

\section{Variational Characterizations}\label{RSWsec:Var}

Delsarte, Goethals and Seidel~\cite{DelGS77} showed that
$X_N = \{ \bsx_1, \ldots, \bsx_N\}\subset \mathbb{S}^{d}$
is a spherical $t$-design if and only if the Weyl sums satisfy
\begin{equation}\label{RSWeq:Wsum}
   r_{\ell,k}(X_N) := \sum_{j=1}^N Y_{\ell,k} (\bsx_j) = 0
 \qquad k=1, \ldots, Z(d,\ell), \quad \ell = 1,\ldots, t ,
\end{equation}
as the integral of all spherical harmonics of degree $\ell \geq 1$ is zero
from orthogonality with the constant ($\ell = 0$) polynomial $Y_{0,1} = 1$
which is not included.

In matrix form
\[
  \bsr(X_N) := \overline{\bsY} \bse = \bszero
\]
where $\bse = (1,\ldots,1)^T\in\mathbb{R}^{N}$ and
$\overline{\bsY} \in \mathbb{R}^{D(d,t)-1 \times N}$
is the spherical harmonic basis matrix excluding the first row.

Let $\psi_t: [-1, 1] \to \mathbb{R}$ be a polynomial of degree $t\geq 1$ with
\begin{equation}\label{RSWeq:psi}
  \psi_t(z) = \sum_{\ell=1}^t a_\ell P^{(d+1)}_{\ell}(z), \qquad
  a_\ell > 0 \mbox{ for } \ell = 1,\ldots,t,
\end{equation}
so the generalised Legendre coefficients $a_\ell$
for degrees $\ell=1,\ldots,t$ are all strictly positive.
Clearly any such function $\psi_t$ can be scaled by an arbitrary positive
constant without changing these properties.

Consider now an arbitrary set $X_N$ of $N$ points on $\mathbb{S}^{d}$.
Sloan and Womersley~\cite{SloWom09} considered the variational form
\[
  V_{t,N,\psi}(X_N) :=
  \frac{1}{N^2} \sum_{i=1}^N \sum_{j=1}^N \psi_t (\bsx_i \cdot \bsx_j)
\]
which from (\ref{RSWeq:LegBnd}) satisfies
\[
  0 \le V_{t,N,\psi}(X_N) \le \sum_{\ell=1}^{t} a_\ell = \psi_t (1).
\]
Moreover the average value is
\[
  \overline{V}_{t, N, \psi} :=
  \int_{\mathbb{S}^{d}} \cdots \int_{\mathbb{S}^{d}}
   V_{t,N,\psi} (\bsx_1, \ldots, \bsx_N) d\sigma_d(\bsx_1) \cdots d\sigma_d(\bsx_N) =
   \frac{\psi_t(1)}{N}.
\]
As the upper bound and average of $V_{t,N,\psi}(X_N)$ depend on
$\psi_t(1)$, we concentrate on functions $\psi$ for which $\psi_t(1)$
does not grow rapidly with $t$.

From the addition theorem~(\ref{RSWeq:addthm}), $V_{t,N,\psi}(X_N)$
is a weighted sum of squares with strictly positive coefficients
\begin{equation}
  V_{t,N,\psi}(X_N) =
  \frac{1}{N^2}
   \sum_{\ell=1}^t \frac{a_\ell}{Z(d,\ell)}
   \sum_{k=1}^{Z(d,\ell)} \left(r_{\ell,k} (X_N)\right)^2
   = \frac{1}{N^2} \bsr(X_N)^T\; \bsD \; \bsr(X_N),
   \label{RSWeq:ssq}
\end{equation}
where $\bsD$ is the diagonal matrix with strictly positive diagonal elements
$\frac{a_\ell}{Z(d,\ell)}$ for $k = 1,\ldots,Z(d,\ell), \ell = 1,\ldots,t$.
Thus, from (\ref{RSWeq:Wsum}), $X_N$ is a spherical $t$-design if and only if
\[
  V_{t,N,\psi}(X_N) = 0.
\]
Moreover, if the \emph{global} minimum of $V_{t,N,\psi}(X_N) > 0$
then there are no spherical $t$-designs on $\mathbb{S}^{d}$ with $N$ points.

Given a polynomial $\widehat\psi_t(z)$ of degree $t$ and
strictly positive Legendre coefficients,
the zero order term may need to be removed to
get $\psi_t(z) = \widehat\psi_t(z) - a_0$ where
for $\mathbb{S}^{d}$ and $\alpha = (d-2)/2$,
\[
  a_0 = \int_{-1}^1 \widehat\psi_t(z) \left(1-z^2\right)^\alpha dz .
\]
Three examples of polynomials on $[-1, 1]$ with strictly positive Legendre coefficients
for $\mathbb{S}^{d}$ and zero constant term, with $\alpha = (d-2)/2$ are:
\begin{itemize}
\item[\textbf{Example 1}]
\begin{equation}\label{RSWeq:psi1}
  \psi_{1,t}(z) = z^{t-1} + z^{t} - a_0
\end{equation}
where
\begin{equation}\label{RSWeq:psi1a0}
  a_0 = 
        \frac{\Gamma(\alpha+3/2)}{\sqrt{\pi}}
        \left\{ \begin{array}{ll}
             \frac{\Gamma(t/2)}{\Gamma(\alpha+1+t/2)} & t \mbox{ odd}, \\[2ex]
             \frac{\Gamma((t+1)/2)}{\Gamma(\alpha+3/2+t/2)} & t \mbox{ even}.
        \end{array} \right.
\end{equation}
For $d = 2$ this simplifies to $a_0 = 1/t$ if $t$ is odd and $a_0 = 1/(t+1)$ if $t$ is even.
This function was used by Grabner and Tichy~\cite{GraTic93} for symmetric  point sets
where only even values of $t$ need to be considered,
as all odd degree polynomials are integrated exactly.
\item[\textbf{Example 2}]
\begin{equation}\label{RSWeq:psi2}
  \psi_{2,t}(z) = \left(\frac{1+z}{2} \right)^t - a_0
\end{equation}
where
\begin{equation}\label{RSWeq:psi2a0}
        a_0 = \frac{2}{\sqrt{\pi}} 4^\alpha \Gamma(\alpha+3/2)
        \frac{\Gamma(\alpha+1+t)}{\Gamma(2\alpha+2+t)}.
\end{equation}
For $d = 2$ this simplifies to $a_0 = 1/(1+t)$.
This is a scaled version of the function $(1+z)^t$ used by Cohn and Kumar~\cite{CohKum07}
for which $a_0$ must be scaled by $2^t$ producing more cancellation errors
for large $t$.
\item[\textbf{Example 3}]
\begin{equation}\label{RSWeq:psi3}
  \psi_{3,t}(z) = P_t^{(\alpha+1,\alpha)}(z) - a_0
\end{equation}
where $a_0$ is given by (\ref{RSWeq:psi2a0}).
The expansion in terms of Jacobi polynomials in Szeg{\H{o}}~\cite[Section 4.5]{Szego75} gives
\[
    \sum_{\ell=0}^t Z(d,\ell) P^{(d+1)}_{\ell}(z) = \frac{1}{a_0} P_t^{(\alpha+1,\alpha)}(z).
\]
For $S^2$ this is equivalent to
\[
  \sum_{\ell=1}^t (2\ell + 1) P^{(d+1)}_{\ell}(z) = (t+1) P_t^{(1,0)}(z) - 1
\]
used in Sloan and Womersley~\cite{SloWom09}.
\end{itemize}

\section{Quadrature Error}\label{RSWsec:CubErr}

The error for numerical integration depends on the smoothness
of the integrand.
Classical results are based on the error of best approximation
of the integrand $f$ by polynomials~\cite{Ragozin71},
(see also \cite{HesSW10} for more details on $\mathbb{S}^{2}$).
For $f\in C^\kappa(\mathbb{S}^{d})$,
there exists a constant $c=c(\kappa,f)$ such that the numerical integration error satisfies
\[
  \left|\int_{\mathbb{S}^{d}} f(\bsx) \mbox d \sigma_d(\bsx) - \frac{1}{N} \sum_{j=1}^N f(\bsx_j)\right|
  \leq c \; t^{-\kappa} .
\]
If $N = O(t^d)$ then the right-hand-side becomes $N^{-\kappa/d}$.
Thus for functions with reasonable smoothness it pays to
increase the degree of precision $t$.

Similar results are presented in~\cite{BraSSW14},
building on the work of \cite{HesSlo06,Hesse06}, for functions $f$
in a Sobolev space $\mathbb{H}^s(\mathbb{S}^{d})$, $s > d/2$.
The \emph{worst-case-error} for equal weight (quasi Monte-Carlo)
numerical integration using an arbitrary point set $X_N$ is
\begin{equation}
  WCE(X_N, s, d) :=
  \sup_{f\in \mathbb{H}^s(\mathbb{S}^{d}), ||f||_{\mathbb{H}^s(\mathbb{S}^{d})\leq 1}} \
  \left| \int_{\mathbb{S}^{d}} f(\bsx) \mbox d \sigma(\bsx) -
        \frac{1}{N} \sum_{j=1}^N f(\bsx_j) \right|.
  \label{RSWeq:WCE}
\end{equation}
From this it immediately follows that the error for numerical integration satisfies
\[
  \left| \int_{\mathbb{S}^{d}} f(\bsx) \mbox d \sigma_d(\bsx) -
        \frac{1}{N} \sum_{j=1}^N f(\bsx_j) \right| \leq
  WCE(X_N, s, d) \  \| f \|_{\mathbb{H}^s(\mathbb{S}^{d})}.
\]
Spherical $t$-designs $X_N$ with $N = O(t^d)$ points satisfy the optimal order
rate of decay of the worst case error, for any $s > d/2$, namely
\[
  WCE(X_N, s, d)  = O\left(N^{-s/d}\right), \qquad N \to \infty.
\]
Thus spherical $t$-designs with $N = O(t^d)$ points are ideally suited to
the numerical integration of smooth functions.

\section{Computational Issues}\label{RSWsec:Comp}

The aim is to find a spherical $t$-design with $N$ points on $\mathbb{S}^{d}$
by finding a point set $X_N$ achieving the global minimum of zero
for the variational function $V_{t, N, \psi}(X_N)$.
This section considers several computational issues:
the evaluation of $V_{t, N, \psi}(X_N)$ either as a double sum
or using its representation (\ref{RSWeq:ssq}) as a sum of squares;
the parametrisation of the point set $X_N$;
the number of points $N$ as a function of $t$ and $d$;
the choice of optimization algorithm which requires evaluation of
derivatives with respect to the chosen parameters;
exploiting the sum of squares structure which requires evaluating the
spherical harmonics and their derivatives;
and imposing structure on the point set, for example symmetric (antipodal) point sets.
An underlying issue is that optimization problems with points on
the sphere typically have many different
local minima with different characteristics.
Here we are seeking both a global minimizer with value $0$ and one with good geometric
properties as measured by the mesh ratio.

The calculations were performed using Matlab, on a Linux computational cluster using nodes
with up to 16 cores.
In all cases analytic expressions for the derivatives with respect to the chosen parametrisation 
were used.

\subsection{Evaluating Criteria}\label{RSWsec:Eval}

Although the variational functions are nonnegative,
there is significant cancellation between the (constant) diagonal elements $\psi_t(1)$
and all the off-diagonal elements with varying signs as
\[
  V_{t,N,\psi}(X_N) =
  \frac{1}{N} \psi_t(1) + \sum_{i=1}^N \sum_{\stackrel{j=1}{j\neq i}}^N \psi_t(\bsx_i\cdot\bsx_j).
\]
\begin{figure}[ht]
\centering
\begin{tabular}{cc}
\includegraphics[scale=0.38]{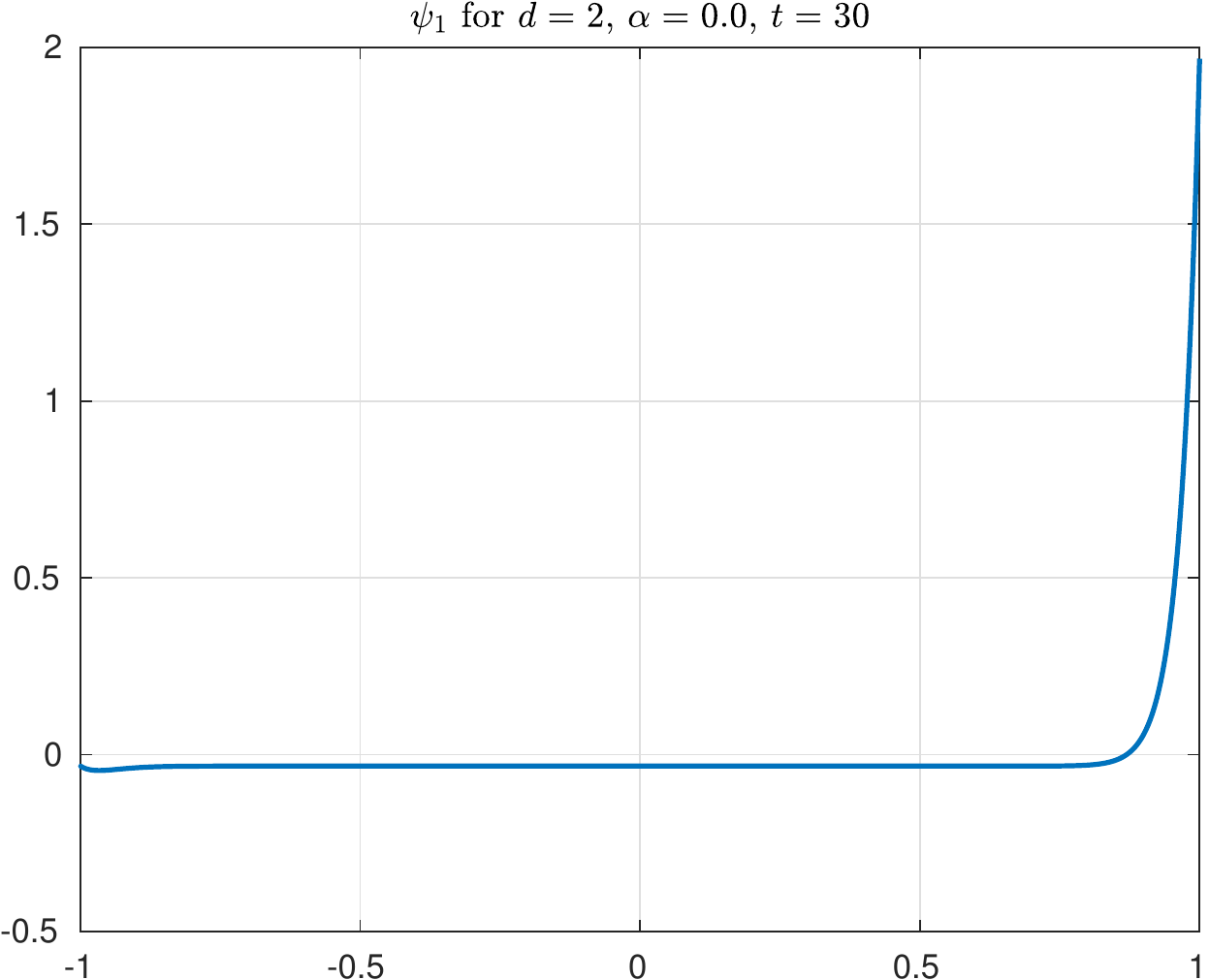} ~~ &
\includegraphics[scale=0.16]{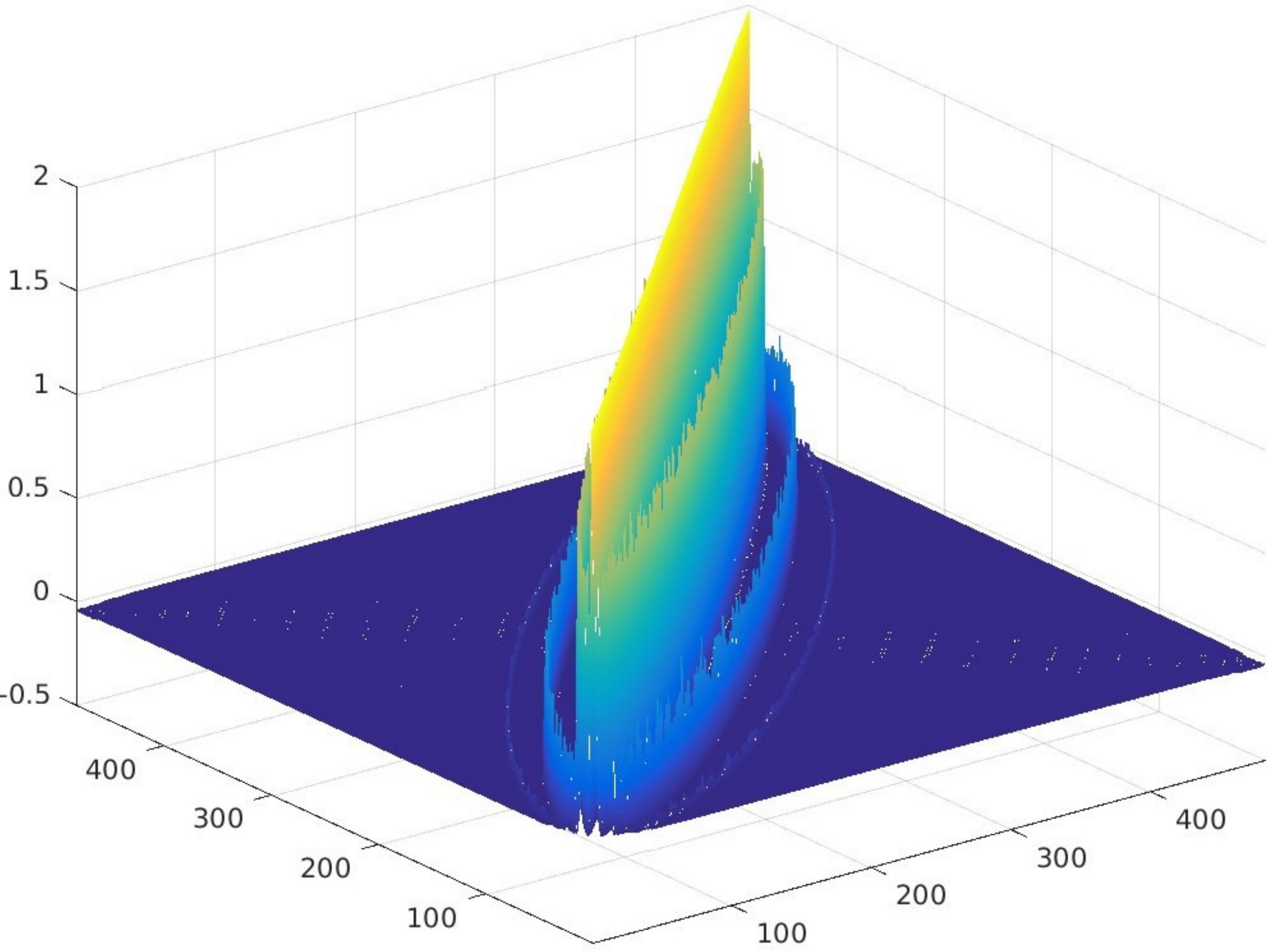}  \\[1.5ex]
\includegraphics[scale=0.38]{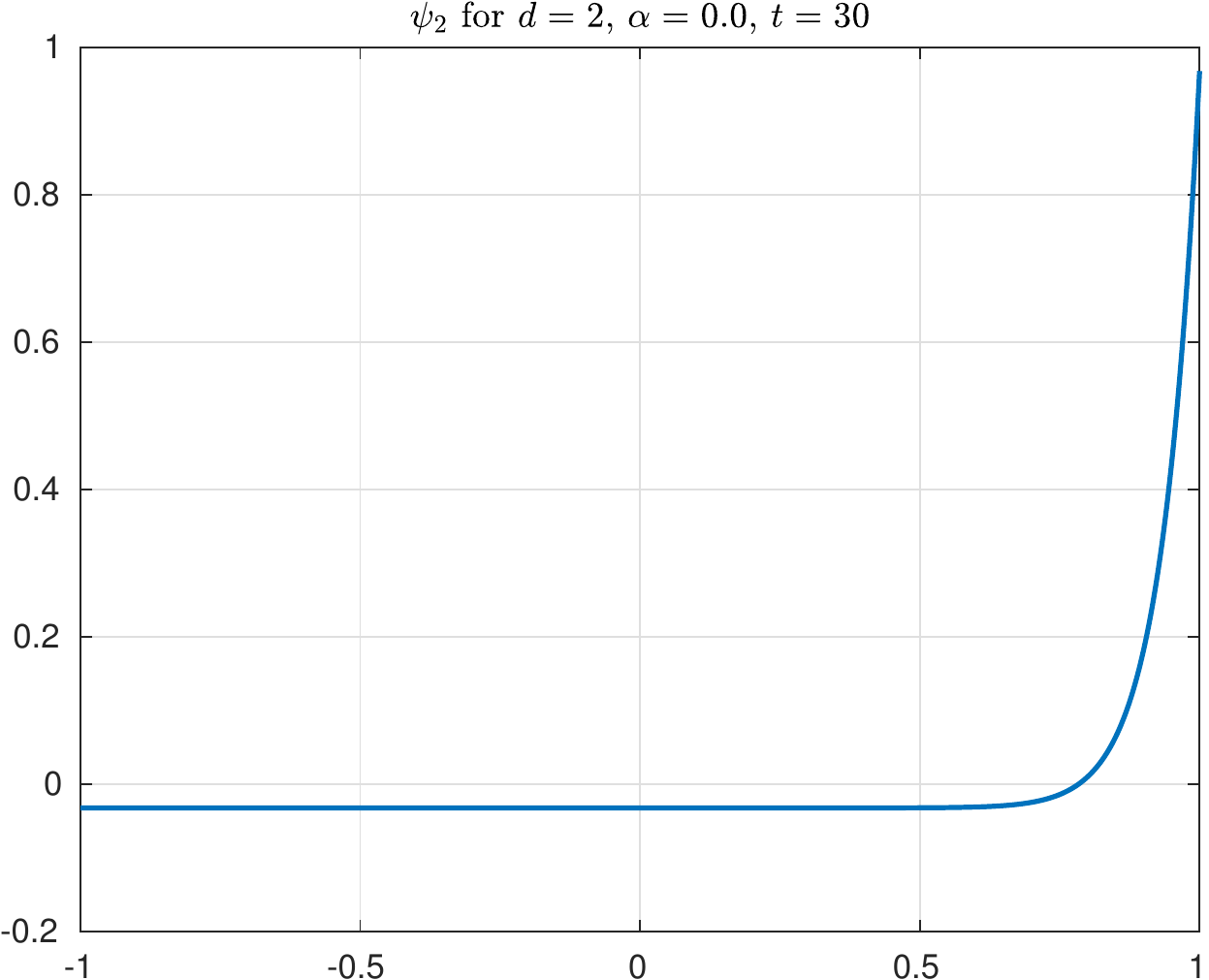} ~~ &
\includegraphics[scale=0.16]{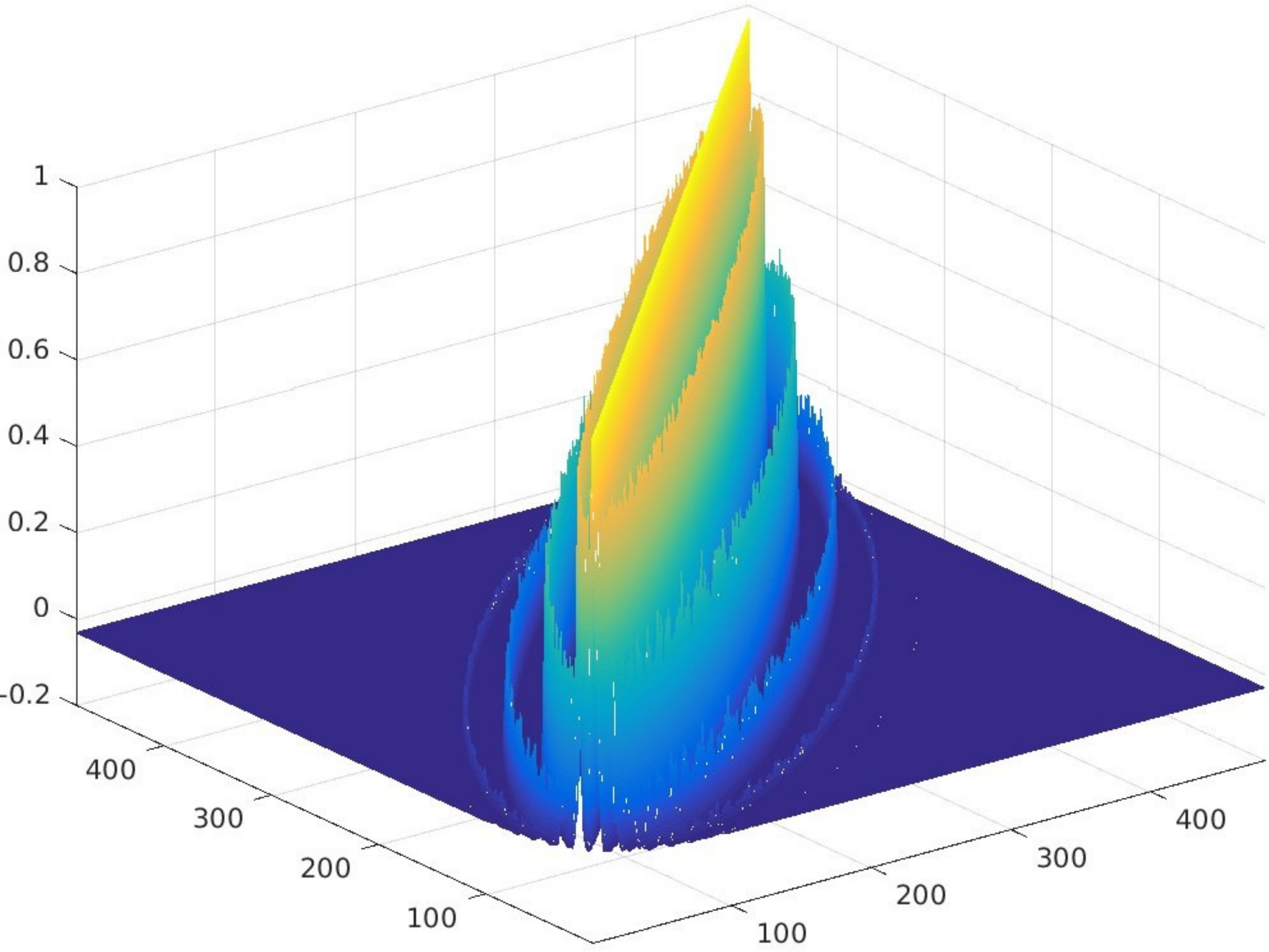}  \\[1.5ex]
\includegraphics[scale=0.38]{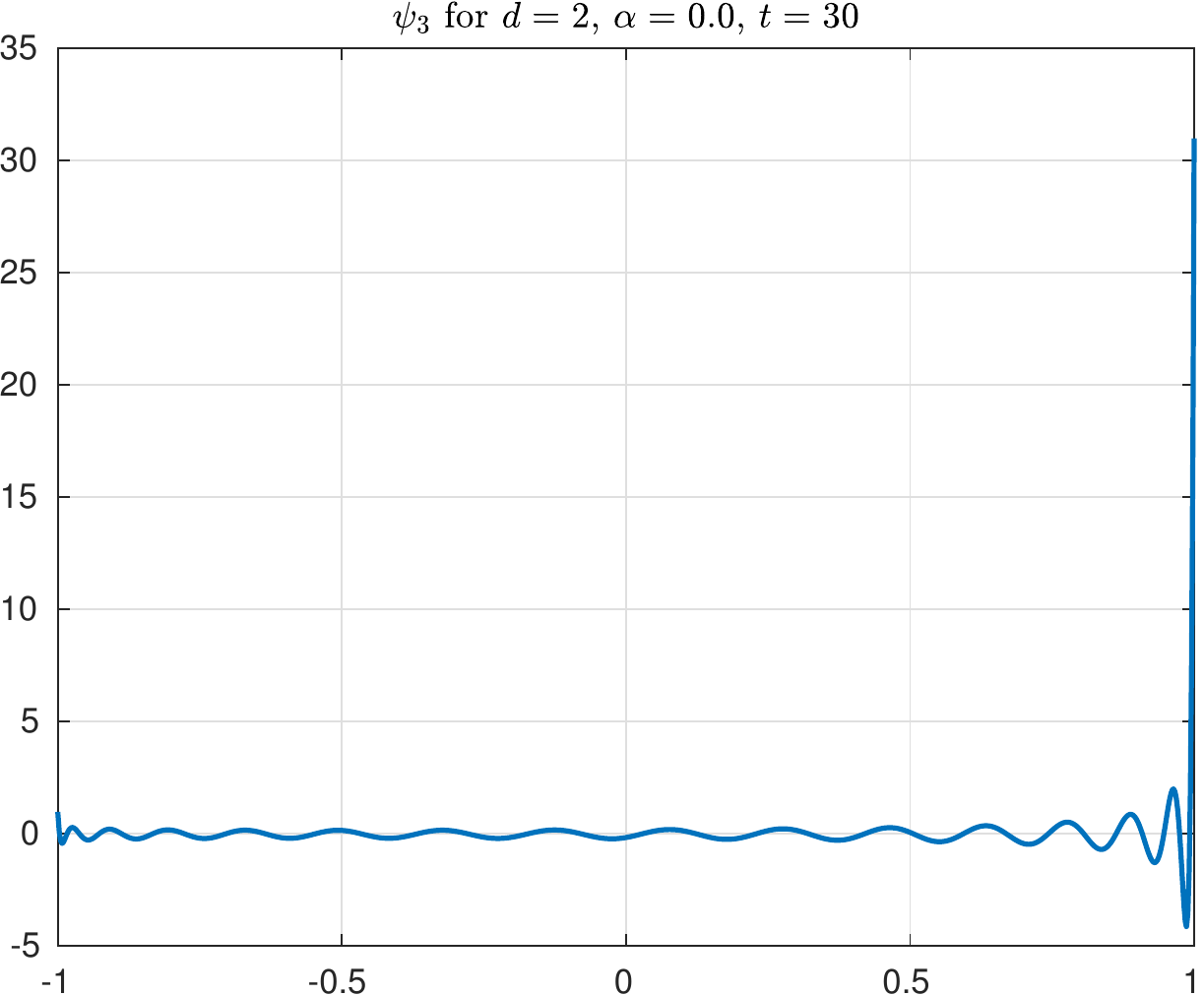} ~~ &
\includegraphics[scale=0.16]{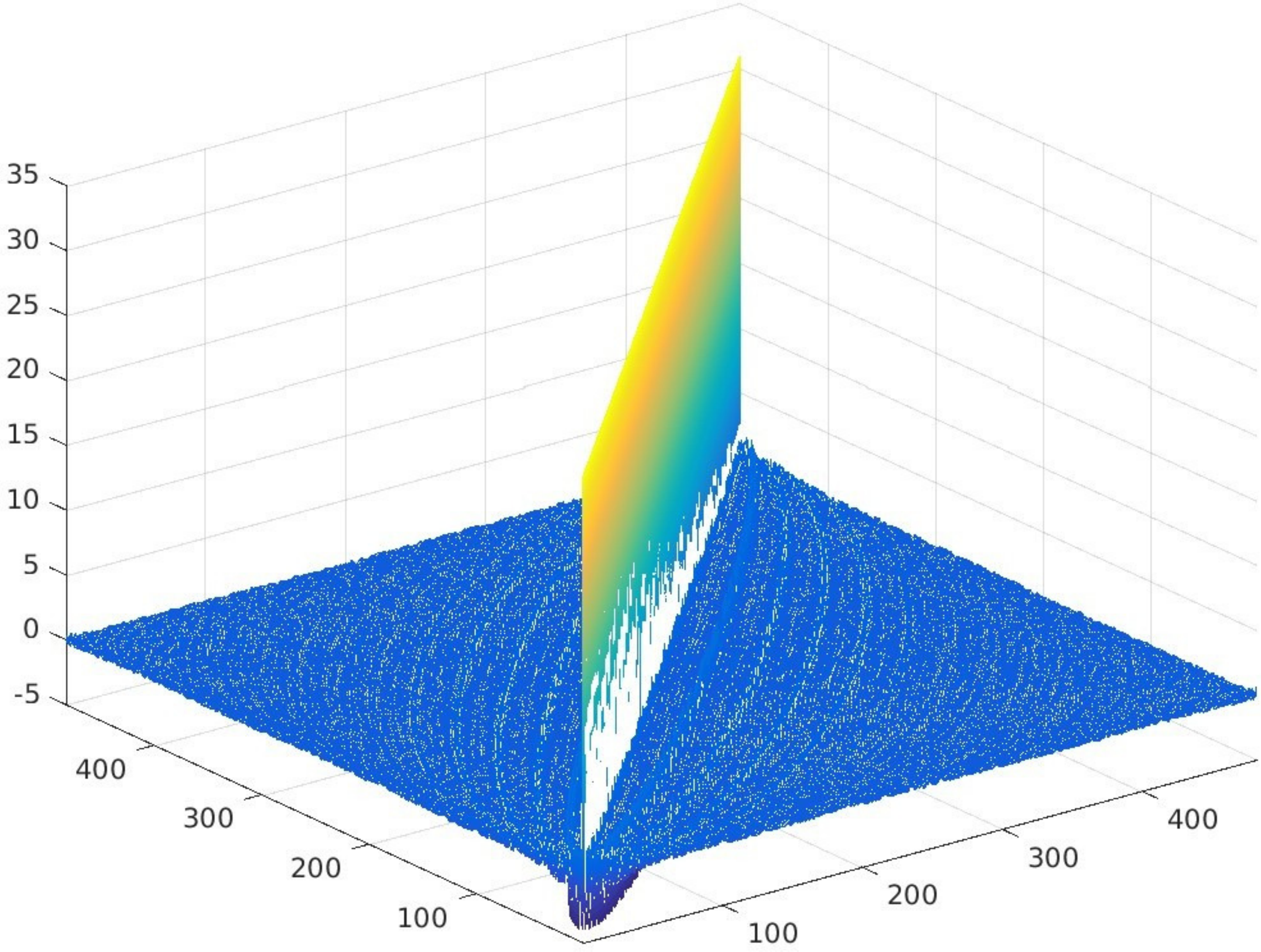}
\end{tabular}
\caption{For $d = 2$, $t = 30$, a spherical $t$-design with $N = 482$,
the functions $\psi_{k,t}$  and arrays $\psi_{k,t}(\bsx_i\cdot\bsx_j)$ for $k = 1, 2, 3$.}
\label{RSWfig:psi}
\end{figure}
Accurate calculation of such sums is difficult, see \cite{Higham96} for example,
especially getting reproducible results on multi-core architecture with dynamic scheduling of parallel
non-associative floating point operations~\cite{DemNgu15}.
Example 1 has $\psi_{1,t}(1) = 2$ and Example 2 has $\psi_{2,t}(1) = 1$,
both independent of $t$, while Example 3 has
\[
  \psi_{3,t}(1) = \frac{\Gamma(t+\alpha+2)}{\Gamma(t+1)\Gamma(\alpha+2)} - 1,
\]
which grows with the degree $t$ (for $d = 2$, $\psi_{3,t}(1) = t$).
These functions are illustrated in Fig.~\ref{RSWfig:psi}.
As the variational objectives can be scaled by an arbitrary positive constant,
you could instead have used $\psi_{3,t} \frac{\Gamma(t+1)\Gamma(\alpha+2)}{\Gamma(t+\alpha+2)}$.
Ratios of gamma functions, as in the expressions for $a_0$, should not be evaluated directly,
but rather simplified for small values of $d$ or evaluated using the log-gamma function.
The derivatives, essential for large scale non-linear optimization algorithms, are readily calculated using
\[
  \nabla_{\bsx_k} V_{t,N,\psi}(X_N) =
  2 \sum_{\stackrel{i=1}{i\neq k}}^N \psi_t'(\bsx_i \cdot \bsx_k) \bsx_i
\]
and the Jacobian of the (normalised) spherical parametrisation (see Section~\ref{RSWsec:SphPar}).

Because of the interest in the use of spherical harmonics for
the representation of the Earth's gravitational field there
has been considerable work, see~\cite{HolFea02,JekLK07} and \cite[Section 7.24.2]{GSL} for example,
on the evaluation of high degree spherical harmonics for $\mathbb{S}^{2}$.
For $(x, y, z)^T \in \mathbb{S}^{2}$ the real spherical harmonics
\cite[Chapter 3, Section 18]{Sansone59}
are usually expressed in terms of the coordinates
$z = \cos(\theta)$ and $\phi$.
In terms of the coordinates $(x,\phi_2) = (\cos(\phi_1), \phi_2)$,
see (\ref{RSWeq:sphpar}) below, they are the $Z(2,\ell) = 2\ell+1$ functions
\begin{eqnarray}
   Y_{\ell,\ell+1-k}(x,\phi_2) & := & \hat{c}_{\ell,k}
   (1-x^2)^{k/2} S_\ell^{(k)}(x) \sin(k\phi_2),
   \quad k=1,\ldots,\ell,\nonumber\\
   Y_{\ell,\ell+1}(x,\phi_2) & := & \hat{c}_{\ell,0} S_\ell^{(0)}(x),
   \label{spherical-harmonics}\\
   Y_{\ell,\ell+1+k}(x,\phi_2) & := & \hat{c}_{\ell,k} (1-x^2)^{k/2}
   S_\ell^{(k)}(x) \cos(k\phi_2),
   \quad k=1,\ldots,\ell. \nonumber
\end{eqnarray}
where $S_\ell^{(k)}(x) = \sqrt{\frac{(\ell-k)!}{(\ell+k)!}} P_\ell^k(x)$
are versions of the Schmidt semi-normalised associated Legendre functions
for which stable three-term recurrences exist for high (about $2700$) degrees and orders.
The normalization constants $\hat{c}_{\ell,0}$, $\hat{c}_{\ell,k}$ are,
for normalised surface measure,
\[
  \hat{c}_{\ell, 0} = \sqrt{ 2\ell+1},  \quad
  \hat{c}_{\ell,k} = \sqrt{2} \sqrt{2\ell+1},
  \quad k = 1,\ldots,\ell,
\]
For $\mathbb{S}^{2}$ these expressions can be used to directly evaluate the
Weyl sums (\ref{RSWeq:Wsum}), and hence their sum of squares, and their derivatives.

\subsection{Spherical Parametrisations}\label{RSWsec:SphPar}

There are many ways to organise a spherical parametrisation of $\mathbb{S}^{d}$.
For $\phi_i \in [0, \pi]$ for $i = 1,\ldots,d-1$ and $\phi_d \in [0, 2\pi)$
define $\bsx\in\mathbb{S}^{d}$ by
\begin{align}
 & x_1   = \cos(\phi_1) \\
 & x_i   = \prod_{k=1}^{i-1} \sin(\phi_k) \cos(\phi_i), \quad i = 2,\ldots,d  \\
 & x_{d+1} = \prod_{k=1}^d \sin(\phi_k)
\label{RSWeq:sphpar}
\end{align}
The inverse transformation used is, for $i = 1,\ldots,d-1$
\begin{align}
  & \phi_i = \left\{
  \begin{array}{ll}
  0  &  \quad\mbox{if  $x_k = 0, \quad k = i,\ldots,d+1$}, \\[1ex]
  \cos^{-1}\left( x_i / \sqrt{\sum_{k=i}^{d+1} x_k^2} \right) & \quad\mbox{otherwise}; \\[1ex]
 \end{array}\right. \\
 & \phi_d = \tan^{-1} \left(x_{d+1}/x_d\right).
 \label{RSWeq:sphinv}
\end{align}
The last component can be calculated using the four quadrant \texttt{atan2} function
and periodicity to get $\phi_d \in [0, 2\pi)$.
Spherical parametrisations introduce potential singularities when
$\phi_i = 0$ or $\phi_i = \pi$ for any $i = 1,\ldots,d-1$.

As all the functions considered are zonal,
they are invariant under an orthogonal transformation (rotation).
Thus the point sets are normalised so that the $d+1$ by $N$ matrix
$\bsX = [\bsx_1 \cdots \bsx_N]$ has
\begin{align*}
  & \bsX_{i,j} = 0 \quad\mbox{for } i = j+1,\ldots,d+1, \quad j = 1,\ldots,\min(d,N) \\
  & \bsX_{i,i} \geq 0 \quad\mbox{for } i  = 1,\ldots,\min(d,N).
\end{align*}
The first normalised point is $\bsx_1 = \bse_1 = (1,0,\ldots,0)\in\mathbb{R}^{d+1}$.
Such a rotation can easily be calculated using the $QR$ factorization of
$\bsX$ combined with sign changes to the rows $Q$.
The corresponding normalised spherical parametrisation has
\[
   \Phi_{i,j} = 0 \quad\mbox{for } i = j,\ldots,d, \quad j = 1,\ldots,\min(d,N),
\]
where the $j$th column of $\Phi$ corresponds to the point $\bsx_j$, $j = 1,\ldots,N$.
The optimisation variables are then $\Phi_{i,j}, i = 1,\ldots,\min(j-1,d), \  j = 2,\ldots, N$,
stored as the vector $\bsphi \in \mathbb{R}^n$ where
\begin{equation}
 n = \left\{\begin{array}{cl}
     \frac{N(N-1)}{2} & \mbox{ for } N \leq d, \\[1ex]
      N d - \frac{d(d+1)}{2} & \mbox{ for } N > d,
     \end{array}\right.
 \label{RSWeq:nphi}
\end{equation}
so
\begin{align*}
& \bsphi_p = \Phi_{i,j}, \quad i = 1,\ldots,\min(j-1,d), \quad j = 2,\ldots,N, \\
 & p = \left\{\begin{array}{cl}
        \frac{\min(j-1,d) (\min(j-1,d)-1)}{2} + i & \mbox{ for } j = 2,\ldots,\min(d,N)\\
        \frac{d(d-1)}{2} + (j-d-1)d + i & \mbox{ for } j = d+1,\ldots,N, \quad N > d.
       \end{array}\right.
\end{align*}
It is far easier to work with a spherical parametrisation with bound constraints
than to impose the quadratic constraints $\bsx_j \cdot \bsx_j = 1, j = 1,\ldots,N$,
especially for large $N$.
As the optimization criteria have the effect of moving the points apart,
the use of the normalised point sets reduces difficulties with singularities
at the boundaries corresponding to $\Phi_{i,j} = 0$ or $\Phi_{i,j} = \pi$, $i = 1,\ldots,d-1$.

For $\mathbb{S}^{2}$, these normalised point sets may be rotated (the variable components re-ordered) using
\[
  Q = \begin{bmatrix}
       0 & 1 & 0 \\ 0 & 0 & 1\\ 1 & 0 & 0
       \end{bmatrix}
\]
to get the commonly~\cite{FliMai99,SloWom04,Womersley_ssd_URL} used normalization with the
first point at the north pole and the second on the prime meridian.

A symmetric (or antipodal) point set ($\bsx \in X_N \iff -\bsx \in X_N$)
must have $N$ even, so can be represented as $\bsX = [\overline{\bsX} \  {-\overline{\bsX}}]$ where the
$d+1$ by $N/2$ array of points $\overline{\bsX}$ is normalised as above.

If only zonal function functions depending just on the inner products $\bsx_i \cdot \bsx_j$ are used
then you could use the variables $\bsZ_{i,j} = \bsx_i \cdot \bsx_j$, so
\[
  \bsZ \in \mathbb{R}^{N\times N}, \quad
  \bsZ^T = \bsZ, \quad \bsZ \succeq 0, \quad
  \mbox{diag}(\bsZ) = \bse, \quad \mbox{rank}(\bsZ) = d+1.
\]
where $\bse = (1,\ldots,1)^T\in\mathbb{R}^N$ and $\bsZ \succeq 0$ indicates $\bsZ$ is positive semi-definite.
The major difficulties with such a parametrisation are the number $N(N-1)/2$ of variables
and the rank condition. Semi-definite programming relaxations (without the rank condition)
have been used to get bounds on problems involving points on the sphere (see, for example, \cite{BacVal2008}).

\subsection{Degrees of Freedom for $\mathbb{S}^{d}$}

Using a normalised spherical parametrisation of $N$ points on
$\mathbb{S}^{d} \subset \mathbb{R}^{d+1}$ there are $n = Nd - d(d+1)/2$ variables
(assuming $N\geq d$).
The number of conditions for a $t$-design is
\[
  m = \sum_{\ell=1}^t Z(d,\ell) = D(d,t) - 1 = Z(d+1,t) - 1.
\]
Using the simple criterion that the
number of variables $n$ is at least the number of conditions $m$,
gives the number of points as
\begin{equation}
  \widehat{N}(d,t) := \ceil{\frac{1}{d}\left(Z(d+1,t) + \frac{d(d+1)}{2} -1\right)} .
  \label{RSWeq:Nhat}
\end{equation}
For $\mathbb{S}^{2}$ there are $n = 2N-3$ variables and $m = (t+1)^2 - 1$ conditions
giving
  \[
      \widehat{N}(2,t) := \ceil{(t+1)^2)/2} + 1.
  \]
Grabner and Sloan~\cite{GraSlo13} obtained separation results
for $N$ point spherical $t$-designs when $N \leq \tau \; 2 N^*$ and $\tau < 1$.
For $d = 2$, $\widehat{N}$ is less than twice the lower bound $N^*$ as
\begin{equation*}\label{E:Ndiff2}
     \widehat{N}(2,t) = 2 N^*(2, t) -  t,
\end{equation*}
but the difference is only a lower order term.
The values for $\widehat{N}(2,t)$, $N^*(2,t)$ and the Yudin lower bound $N^+(2,t)$
are available in Tables~\ref{RSWtab:SF1} --  \ref{RSWtab:S3SDS}.

The idea of exploiting symmetry to reduce the number of conditions
that a quadrature rule should satisfy at least goes back
to Sobolev~\cite{Sobolev62}.
For a symmetric point set (both $\bsx_j, -\bsx_j$ in $X_N$)
then all odd degree polynomials $Y_{\ell,k}$ or $P_\ell^{(d+1)}$
are automatically integrated exactly by an equal weight quadrature rule.
Thus, for $t$ odd, the number of conditions to be satisfied is
\begin{equation}
  m = \sum_{\ell=1}^{(t-1)/2} Z(d,2\ell) = \frac{\Gamma(t+d)}{\Gamma(d+1)\Gamma(t)} - 1.
  \label{RSWeq:m-sym}
\end{equation}
The number of free variables in a normalised symmetric point set
$\bsX = [\overline{\bsX} \quad \mbox{$-\overline{\bsX}$}]$ (assuming $N/2 \geq d$) is
\begin{equation}
  n  = \left( \frac{N d}{2} - \frac{d(d+1)}{2} \right) .
  \label{RSWeq:n-sym}
\end{equation}
Again the simple requirement that $n \geq m$ gives the number of points as
\begin{equation}
  \overline{N}(d,t) := 2
  \ceil{ \frac{1}{d} \left(\frac{\Gamma(t+d)}{\Gamma(d+1)\Gamma(t)} - 1 + \frac{d(d+1)}{2}\right)} .
  \label{RSWeq:Nbar}
\end{equation}
For $d = 2$ this simplifies, again for $t$ odd, to
\[
  \overline{N}(2,t) := 2 \ceil{\frac{t^2 + t + 4}{4}} . 
\]
$\overline{N}(2,t)$ is slightly less than $\widehat{N}(2,t)$, comparable to twice the lower bound $N^*(2,t)$ as
\[
\overline{N}(2,t) = 2N^*(2,t) -\mbox{$\frac{3}{2}$} t +
\left\{ \begin{array}{ll}
       \frac{3}{2}  & \mbox{ if } \mod(t,4) = 1,\\[0.5ex]
       \frac{1}{2}  & \mbox{ if } \mod(t,4) = 3.
       \end{array}\right.
\]
However $\overline{N}(2,t)$ is not less than $\tau \; 2 N^*(2,t)$, $\tau < 1$,
as required by Grabner and Sloan~\cite{GraSlo13}. 

The leading term of both $\widehat{N}(d,t)$ and $\overline{N}(d,t)$ is $D(d,t)/d$,
see Table~\ref{RSWtab:Ndt},
where $D(d,t)$ defined in (\ref{RSWeq:dimP}) is the dimension of $\mathbb{P}_t(\mathbb{S}^{d})$.
From (\ref{RSWeq:eff}), a spherical $t$-design with $\widehat{N}(d,t)$
or $\overline{N}(d,t)$ points has efficiency $E \approx 1$.
Also the leading term of both $\overline{N}(d,t)$ and  $\widehat{N}(d,t)$
is  $2^d/d$ times the leading term of the lower bound $N^*(d,t)$.
\begin{table}[ht]
\[
\begin{array}{l||c|c|c|c}
d & N^*(d,t) & \overline{N}(d,t) & \widehat{N}(d,t) & D(d,t) \\[1ex]
\hline
2 & \frac{t^2}{4} + t + O(1) &  \frac{t^2}{2} + \frac{t}{2} + O(1) &
    \frac{t^2}{2} + t + O(1) & t^2+2t+1 \\[1ex]
3 & \frac{t^3}{24} + \frac{3 t^2}{8} + O(t) &
     \frac{t^3}{9} + \frac{t^2}{3} + O(t) & \frac{t^3}{9} + \frac{t^2}{2} + O(t) &
    \frac{t^3}{3} + \frac{3 t^2}{2} + O(t)  \\[1ex]
4 & \frac{t^4}{192} + \frac{t^3}{12} + O(t^2) &
    \frac{t^4}{48} + \frac{t^3}{8} + O(t^2) & \frac{t^4}{48} + \frac{t^3}{6} + O(t^2) & 
    \frac{t^4}{12} + \frac{2 t^3}{3} + O(t^2)\\[1ex]
5 & \frac{t^5}{1920} + \frac{5 t^4}{384} + O(t^3)&
    \frac{t^5}{300} + \frac{t^4}{30} + O(t^3) & \frac{t^5}{300} + \frac{t^4}{24} + O(t^3) &
    \frac{t^5}{60} + \frac{5 t^4}{24} + O(t^3)
\end{array}
\]
\caption{The lower bound $N^*(d,t)$, the number of points $\overline{N}(d,t)$ (symmetric point set)
and $\widehat{N}(d,T)$ to match the number of conditions and the dimension of $\mathbb{P}_t(\mathbb{S}^{d})$
for $d = 2, 3, 4, 5$}
\label{RSWtab:Ndt}
\end{table}

\subsection{Optimization Algorithms}\label{RSWsec:Opt}

As  with many optimization problems on the sphere there are many
distinct (not related by an orthogonal transformation or permutation)
points sets giving local minima of the optimization objective.
For example, Erber and Hockney~\cite{ErbHoc97} and Calef et al~\cite{CalGS15}
studied the minimal energy problem for the sphere and
the large number of stable configurations.

Gr{\"a}f and Potts~\cite{GraPot11} develop optimization methods on general Riemannian manifolds,
in particular $\mathbb{S}^{2}$, and both Newton-like and conjugate gradients methods.
Using a fast method for spherical Fourier coefficients at non-equidistant points
they obtain approximate spherical designs for high degrees.

While mathematically it is straight forward to conclude that if
$V_{t,N,\psi}(X_N) = 0$ then $X_N$ is a spherical $t$-design,
deciding when a quantity is zero with the limits of standard double precision
floating point arithmetic with machine precision $\epsilon = 2.2\times 10^{16}$
is less clear (should $10^{-14}$ be regarded as zero?).
Extended precision libraries and packages like Maple or Mathematica can help.
A point set $X_N$ with  $V_{t,N,\psi}(X_N) \approx\epsilon$
does not give a mathematical proof that is $X_N$ is a spherical $t$-design,
but $X_N$ may still be computationally useful in applications.

On the other hand showing that the global minimum if $V_{t,N,\psi}(X_N)$
is strictly positive, so no spherical $t$-design with $N$ points exist,
is an intrinsically hard problem problem.
Semi-definite programming~\cite{VanBoy96} provides an approach~\cite{ParStu03}
to the global optimization of polynomial sum of squares for modest degrees.

For $d = 2$ a variety of gradient based bound constrained optimization methods,
for example the limited memory algorithm \cite{ByrLNZ95,MorNoc11},
were tried both to minimise the variational forms $V_{t,N,\psi}(X_N)$.
Classically, see \cite{NocWri06} for example, methods can exploit the sum of squares
structure  $\bsr(X_N)^T \bsr(X_N)$.
In both cases it is important to provide derivatives of the objective
with respect to the parameters.
Using the normalised spherical parametrisation $\bsphi$ of $X_N$,
the Jacobian of the residual $\bsr(\bsphi)$ is $\bsA: \mathbb{R}^{n} \to \mathbb{R}^{m \times n}$
where $n = dN - d(d+1)/2$ and $m = D(d,t)-1$
\[
   \bsA_{i,j}(\bsphi) = \frac{\partial r_i(\bsphi)}{\partial \phi_j}, 
   \qquad i = 1,\ldots,m, \quad j = 1,\ldots,n,
\]
where $i = (\ell-1) Z(d+1,\ell-1) + k$, for
$k = 1,\ldots,Z(d,\ell), \ell =1,\ldots,t$.

For symmetric point sets with $N = \overline{N}(d,t)$ points,
the number of variables $n$ is given by (\ref{RSWeq:n-sym})
and the number of conditions $m$ by (\ref{RSWeq:m-sym})
corresponding to even degree spherical harmonics.

The well-known structure of a nonlinear least squares problem,
see \cite{NocWri06} for example, gives, ignoring the $1/N^2$ scaling in (\ref{RSWeq:ssq}),
\begin{align}
  & f(\bsphi) = \bsr(\bsphi)^T \; \bsD \; \bsr(\bsphi), \\
  & \nabla f(\bsphi) = 2 \bsA(\bsphi)^T \bsD \, \bsr(\bsphi), \\
  & \nabla^2 f(\bsphi) = 2 \bsA(\bsphi)^T \bsD \bsA(\bsphi) + 
    2 \sum_{i=1}^m r_i(\bsphi) D_{ii} \nabla^2 r_i(\bsphi).
  \label{RWeq:ssq}
\end{align}
If $\bsphi^*$ has $\bsr(\bsphi^*) = \bszero$ and $\bsA(\bsphi^*)$ has rank $n$,
the Hessian $\nabla^2 f(\bsphi^*) = 2 \bsA(\bsphi^*)^T \bsD \bsA(\bsphi^*)$ is positive definite
and $\bsphi^*$ is a strict global minimizer.
Here this is only possible when $n = m$,
for example when $d = 2$ and $t$ is odd,
see Tables~\ref{RSWtab:SF1}, \ref{RSWtab:SF2}, \ref{RSWtab:SF3},
and in the symmetric case when $t \mod 4 = 3$, 
see Tables~\ref{RSWtab:SS1}, \ref{RSWtab:SS2}, \ref{RSWtab:SS3}.
For $d = 2$ the other values of $t$ have $n = m + 1$, so there is generically 
a one parameter family of solutions even when the Jacobian has full rank.
When $d = 3$, the choice $N = \widehat{N}(3,t)$ gives $n = m$, $n = m + 1$ or $n = m+3$
depending on the value of $t$, see Table~\ref{RSWtab:S3SDC}.
Thus a Levenberg-Marquadt or trust region method, see \cite{NocWri06} for example,
in which the search direction satisfies
\[
   \left(\bsA^T \bsD \bsA + \nu \bsI \right) \bsd = \bsA^T \bsD \bsr
\]
was used. 
When $n > m$ the Hessian of the variational form $V_{t,N,\psi}(X_N)$ evaluated using one
of the three example functions (\ref{RSWeq:psi1}), (\ref{RSWeq:psi2}) or (\ref{RSWeq:psi3})
will also be singular at the solution.
These disadvantages could have been reduced by choosing the number of points $N$
so that $n < m$, but then there may not be solutions with $V_{t,N,\psi}(X_N) = 0$,
that is spherical $t$-designs may not exist for that number of points.

Many local solutions were found as well as (computationally) global solutions
which differed depending on the starting point and the algorithm parameters
(for example the initial Levenberg-Marquadt parameter $\nu$, initial trust region,
line search parameters etc).
Even when $n = m$ there are often multiple spherical designs for same $t$, $N$,
which are strict global minimisers, but have different
inner product sets $\mathcal{A}(X_N)$ in (\ref{RSWeq:IPset}) and different mesh ratios.

\subsection{Structure of Point Sets}\label{RSWsec:Struc}

There are a number of issues with the spherical designs studied here.
\begin{itemize}
\item
There is no proof that spherical $t$-designs on $\mathbb{S}^{d}$ with $N = t^d/d + O(t^{d-1})$
points exist for all $t$ (that is the constant in the Bondarenko et al result~\cite{BonRV13}
is $C_d = 1/d$ (or lower), as suggested by \cite{HarSlo96} for $\mathbb{S}^{2}$).
\item
The point sets are not nested, that is the points of a spherical $t$-design
are not necessarily a subset of the points of a $t'$-design for some $t' > t$.
\item
The point sets do not lie on bands of equal $\phi_1$ (latitude on $\mathbb{S}^{2}$)
making them less amenable for FFT based methods.
\item
The point sets are obtained by extensive calculation, rather than generated by a simple
algorithms as for generalized spiral or equal area points on $\mathbb{S}^{2}$~\cite{RakSZ94}.
Once calculated the point sets are easy to use.
\end{itemize}
An example of a point set on $\mathbb{S}^{2}$ that satisfies the last three issues are
the HELAPix points\cite{HEALPix05}, which provide a hierarchical, equal area (so exact for constants),
iso-latitude set of points widely used in cosmology.

\appendix

\section{Tables of Results}\label{RSWsec:Tables}

\subsection{Spherical $t$-Designs with no Imposed Symmetry for $\mathbb{S}^{2}$}

\begin{figure}[ht]
\centering
\includegraphics[scale=0.70]{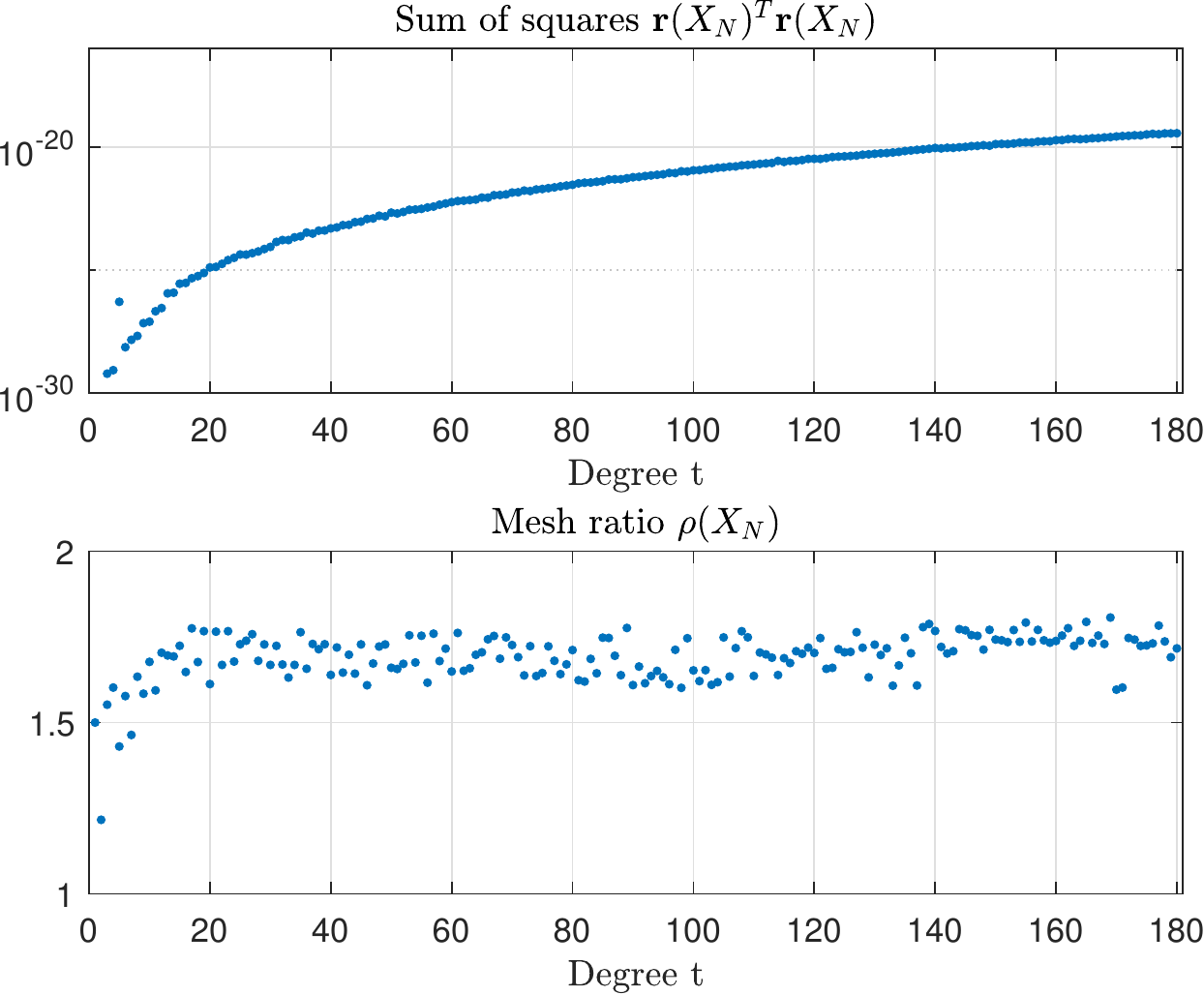}
\caption{Sum of squares of Weyl sums and mesh ratios for spherical $t$-designs on $\mathbb{S}^{2}$}
\label{RSWfig:SF}
\end{figure}
From Tables~\ref{RSWtab:SF1}, \ref{RSWtab:SF2} and \ref{RSWtab:SF3} the variational criteria
based on the three functions $\psi_{1,t}$, $\psi_{2,t}$ and $\psi_{3,t}$
all have values close to the double
precision machine precision of $\epsilon = 2.2\times 10^{-16}$ for all degrees $t = 1,\ldots,180$.
Despite being theoretically non-negative, rounding error sometimes gives negative
values, but still close to machine precision.
The potential values using $\psi_{3,t}$ are slightly larger due to the larger value of $\psi_{3,t}(1)$.
The tables also give the unscaled sum of squares $\bsr(X_N)^T\bsr(X_N)$,
which is also plotted in Fig.~\ref{RSWfig:SF}.
These tables also list both the Delsarte, Goethals and Seidel lower bounds $N^*(2,t)$
and the Yudin lower bound $N^+(2,t)$, plus the actual number of points $N$.
The number of points $N = \widehat{N}(2,t)$,
apart from  $t = 3, 5, 7, 9, 11, 13, 15$ when $N = \widehat{N}(2,t)-1$.
There may well be spherical $t$-designs with smaller values of $N$ and
special symmetries, see \cite{HarSlo96} for example.
For all these point sets the mesh ratios $\rho(X_N)$ are less than $1.81$,
see Fig.~\ref{RSWfig:SF}.
All these point sets are available from~\cite{Womersley_ssd_URL}.

\subsection{Symmetric Spherical $t$-Designs for $\mathbb{S}^{2}$}

\begin{figure}[ht]
\centering
\includegraphics[scale=0.70]{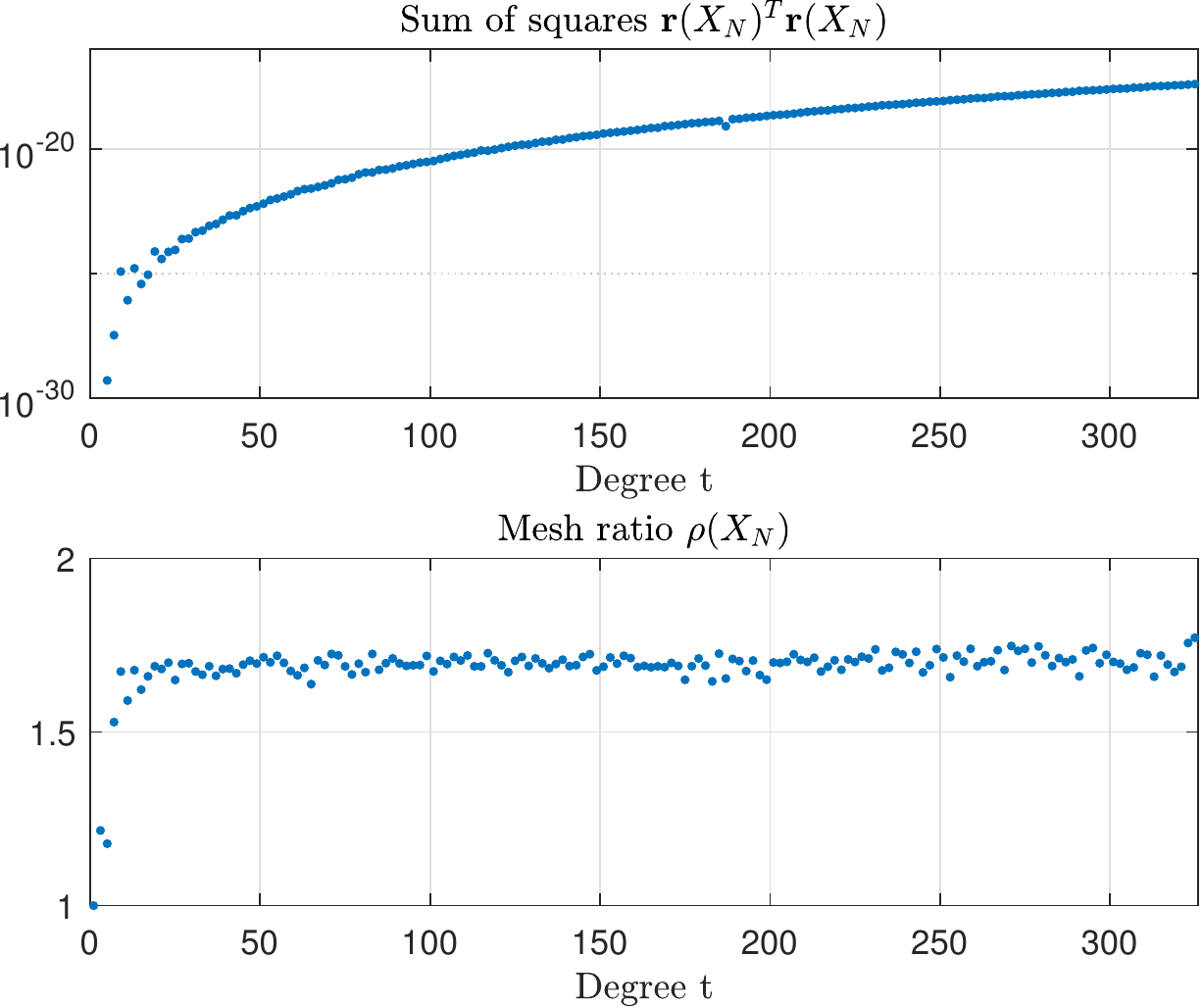}
\caption{Sum of squares of Weyl sums and mesh ratios for symmetric spherical $t$-designs on $\mathbb{S}^{2}$}
\label{RSWfig:SS}
\end{figure}
For $\mathbb{S}^{2}$ a $t$-design with a sightly smaller number of points $\overline{N}(2,t)$
can be found by constraining the point sets to be symmetric (antipodal).
A major computational advantage of working with symmetric point sets
is the reduction (approximately half), for a given degree $t$,
in the number of optimization variables $n$ and the number of terms $m$ in the Weyl sums.
Tables ~\ref{RSWtab:SS1}, \ref{RSWtab:SS2} and \ref{RSWtab:SS3}
list the characteristics of the calculated $t$-designs for $t = 1,3,5,\ldots,325$,
as a symmetric $2k$-design is automatically a $2k+1$-design.
These tables have $t = \overline{N}(2,t)$ except for $t = 1, 7, 11$.
These point sets, again available from~\cite{Womersley_ssd_URL},
provide excellent sets of points for numerical integration on $\mathbb{S}^{2}$
with mesh ratios all less than $1.78$ for degrees up to $325$,
as illustrated in Fig.~\ref{RSWfig:SS}.

\subsection{Designs for $d =  3$}

For $d = 3$, $Z(3,\ell) = (\ell+1)^2$, so the dimension of the space
of polynomials of degree at most $t$ in $\mathbb{S}^{3}$ is
$D(3,t) = Z(4,t) = (t+1)(t+2)(2t+3)/6$.
Comparing the number of variables with the number of conditions,
with no symmetry restrictions, gives
\[
  \widehat{N}(3,t) = \ceil{\frac{2 t^3 + 9 t^2 + 13 t + 36}{18}},
\]
while for symmetric spherical designs on $\mathbb{S}^{3}$
\[
  \overline{N}(3,t) = 2 \ceil{\frac{t^3 + 3 t^2 + 2 t + 30}{18}}.
\]

\begin{figure}[ht]
\centering
\includegraphics[scale=0.40]{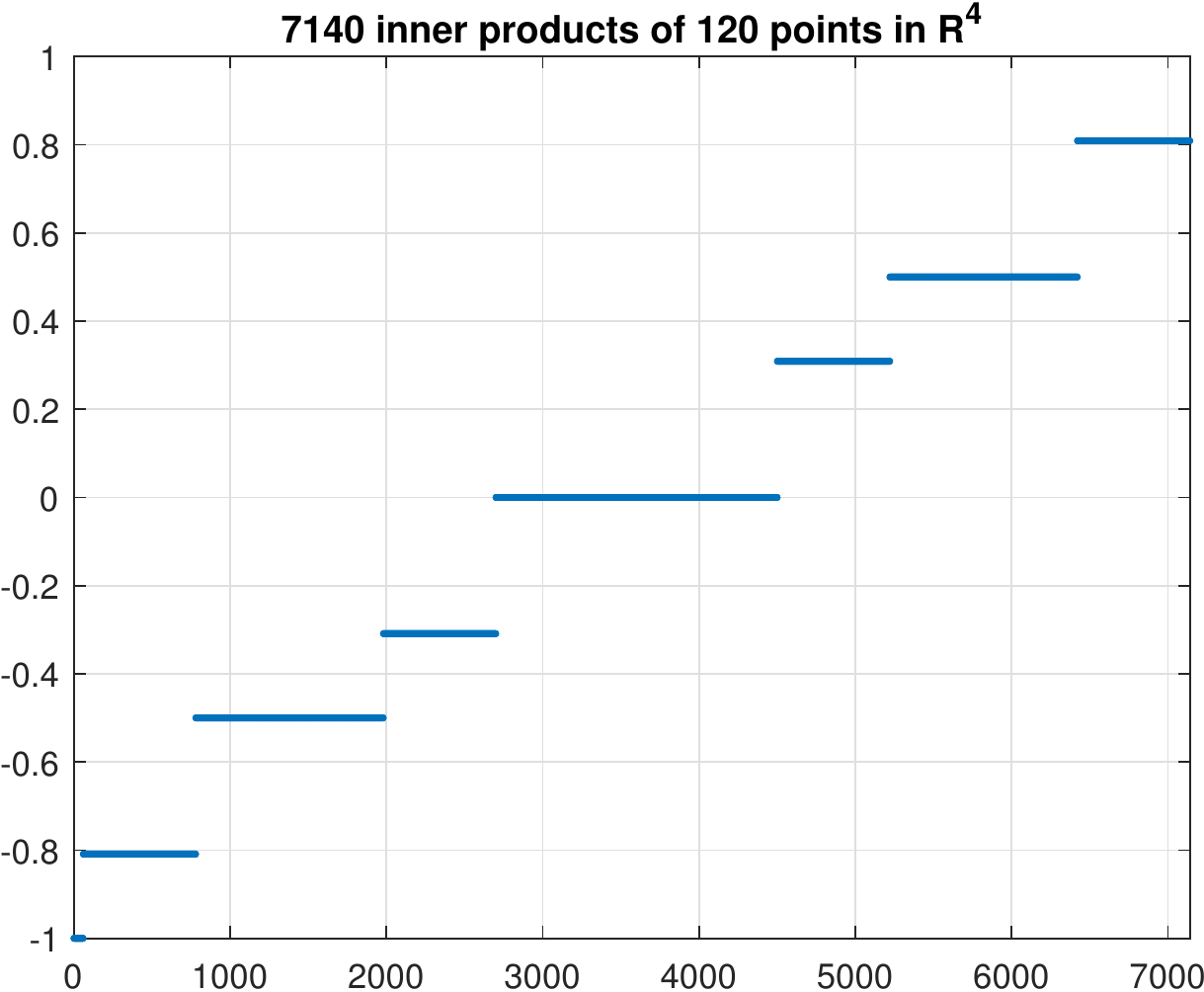}\qquad 
\includegraphics[scale=0.40]{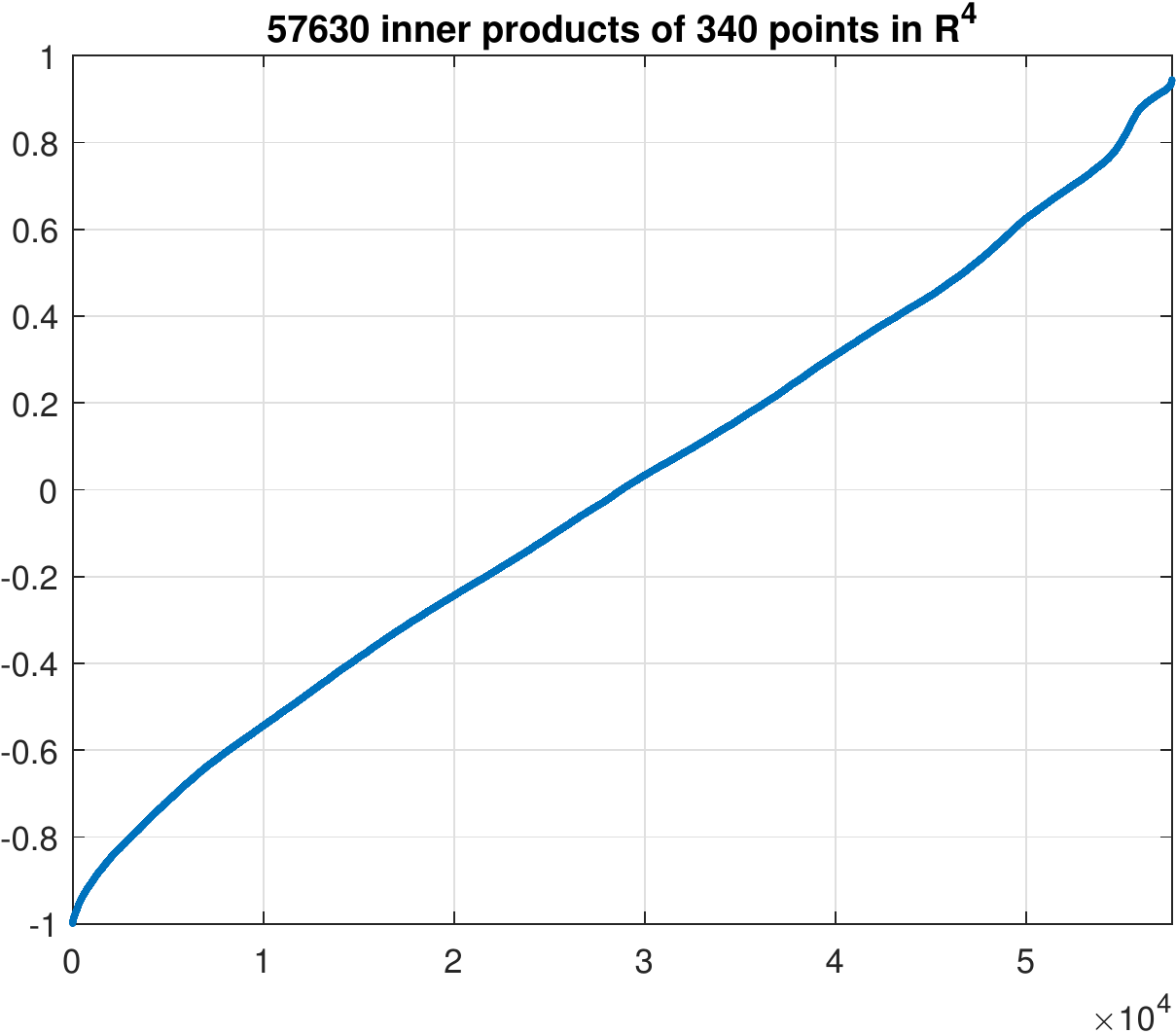}
\caption{Inner product sets $\mathcal{A}(X_N)$ for
$600$-cell with $N = 120$  and $13$-design with $N = 340$ on $\mathbb{S}^{3}$}
\label{RSWfig:S3ip}
\end{figure}

The are six regular convex polytopes with $N = 5, 8, 16, 24, 120$ and $600$
vertices on $S^3$~\cite{Coxeter73}
(the $5$-cell, $16$-cell, $8$-cell, $24$-cell, $600$-cell and $120$-cell respectively)
giving spherical $t$-designs for for $t = 2, 3, 5, 7, 9, 11$ and $11$.
The energy of regular sets on $\mathbb{S}^{3}$ with $N = 2, 3, 4, 5, 6 ,8, 10, 12, 13, 24, 48$
has been studied by \cite{AgbKGL15}.
The $N = 24$ vertices of the D4 root system~\cite{CohCEK07} provides a one-parameter family
of $5$-designs on $\mathbb{S}^{3}$.
The Cartesian coordinates of the regular point sets are known,
and these can be numerically verified to be spherical designs.
The three variational criteria using (\ref{RSWeq:psi1}), (\ref{RSWeq:psi2}) and (\ref{RSWeq:psi3})
are given for these point sets in Table~\ref{RSWtab:S3SDR}.
Fig.~\ref{RSWfig:S3ip} clearly illustrates the difference between
the widely studied~\cite{DelGS77,ConSlo99,BoyDel15}
inner product set $\mathcal{A}(X_N$) for a regular point set (the $600$-cell with $N = 120$)
and a computed spherical $13$-design with $N = 340$.  

The results of some initial  experiments in minimising the three variational criteria 
are given in Tables~\ref{RSWtab:S3SDC} and \ref{RSWtab:S3SDS}.
For $d > 2$, it is more difficult to quickly generate a point set with a good mesh ratio
to serve as an initial point for the optimization algorithms.
One strategy is to randomly generate starting points, but this both makes
the optimization problem harder and tends to produce nearby point sets which
are local minimisers and have poor mesh ratios as the random initial points may have small separation~\cite{BraRSSWW16}.
Another possibility is the generalisation of equal area points to $d > 2$
by Leopardi~\cite{Leopardi06}.
For a given $t$ and $N$ there are still many different point
sets with objective values close to $0$ and different mesh ratios.
To fully explore spherical $t$-designs for $d > 2$, a stable implementation
of the spherical harmonics is needed, so that least squares minimisation
can be fully utilised.

\begin{table}[ht]
\begin{equation*}
\begin{array}{|c|ccccc|rrrc|ccc|}
\hline
 t & N^*(2,t) & N^+(2,t) & N & n & m & 
 \scriptstyle V_{t,N,\psi_1}(X_N) & \scriptstyle V_{t,N,\psi_2}(X_N) & \scriptstyle V_{t,N,\psi_3}(X_N) & \scriptstyle \bsr(X_N)^T\bsr(X_N) & 
 \delta(X_N) & h(X_N) & \rho(X_N)\\ 
\hline
  1 &     2 &     2 &     3 &     3 &     3 &
\mbox{3.7e-17} & \mbox{1.9e-17} & \mbox{1.2e-17} & \mbox{3.0e-32} &
\mbox{ 2.0944} & \mbox{ 1.5708} & 1.50
\\ 
  2 &     4 &     4 &     6 &     9 &     8 &
\mbox{3.5e-17} & \mbox{3.5e-17} & \mbox{4.3e-17} & \mbox{1.1e-31} &
\mbox{ 1.5708} & \mbox{ 0.9553} & 1.22
\\ 
  3 &     6 &     6 &     8 &    13 &    15 &
\mbox{3.1e-17} & \mbox{-4.8e-18} & \mbox{-1.0e-17} & \mbox{6.2e-30} &
\mbox{ 1.2310} & \mbox{ 0.9553} & 1.55
\\ 
  4 &     9 &     9 &    14 &    25 &    24 &
\mbox{-2.1e-18} & \mbox{-5.5e-18} & \mbox{4.5e-17} & \mbox{8.6e-30} &
\mbox{ 0.8630} & \mbox{ 0.6913} & 1.60
\\ 
  5 &    12 &    12 &    18 &    33 &    35 &
\mbox{-3.4e-17} & \mbox{6.9e-19} & \mbox{-1.1e-16} & \mbox{5.2e-27} &
\mbox{ 0.8039} & \mbox{ 0.5749} & 1.43
\\ 
  6 &    16 &    16 &    26 &    49 &    48 &
\mbox{2.0e-17} & \mbox{5.0e-18} & \mbox{6.0e-17} & \mbox{7.4e-29} &
\mbox{ 0.6227} & \mbox{ 0.4911} & 1.58
\\ 
  7 &    20 &    20 &    32 &    61 &    63 &
\mbox{1.1e-17} & \mbox{4.5e-18} & \mbox{6.1e-17} & \mbox{1.5e-28} &
\mbox{ 0.5953} & \mbox{ 0.4357} & 1.46
\\
\hline
\end{array}
\end{equation*}
\caption{Spherical $t$-designs on $\mathbb{S}^2$ with no symmetry restrictions,
$N = \widehat{N}(2,t)$ and degrees $t = 1$ -- $60$,
except $t = 3, 5, 7, 9, 11, 13, 15$ when $N = \widehat{N}(2,t)-1$}
\label{RSWtab:SF1}
\end{table}
\begin{table}[ht]
\begin{equation*}
\begin{array}{|c|ccccc|rrrc|ccc|}
\hline
 t & N^*(2,t) & N^+(2,t) & N & n & m & 
 \scriptstyle V_{t,N,\psi_1}(X_N) & \scriptstyle V_{t,N,\psi_2}(X_N) & \scriptstyle V_{t,N,\psi_3}(X_N) &
 \scriptstyle \bsr(X_N)^T\bsr(X_N) & \delta(X_N) & h(X_N) & \rho(X_N)\\ 
\hline
 61 &   992 &  1065 &  1923 &  3843 &  3843 &
\mbox{-9.4e-17} & \mbox{5.3e-17} & \mbox{-4.3e-17} & \mbox{6.4e-23} &
\mbox{ 0.0663} & \mbox{ 0.0584} & 1.76
\\ 
 62 &  1024 &  1099 &  1986 &  3969 &  3968 &
\mbox{2.6e-17} & \mbox{2.4e-17} & \mbox{-1.9e-16} & \mbox{6.7e-23} &
\mbox{ 0.0699} & \mbox{ 0.0577} & 1.65
\\ 
 63 &  1056 &  1134 &  2049 &  4095 &  4095 &
\mbox{2.6e-17} & \mbox{-4.8e-18} & \mbox{-1.5e-16} & \mbox{7.1e-23} &
\mbox{ 0.0680} & \mbox{ 0.0564} & 1.66
\\ 
 64 &  1089 &  1170 &  2114 &  4225 &  4224 &
\mbox{3.3e-17} & \mbox{2.7e-17} & \mbox{6.5e-17} & \mbox{7.4e-23} &
\mbox{ 0.0662} & \mbox{ 0.0562} & 1.70
\\ 
 65 &  1122 &  1206 &  2179 &  4355 &  4355 &
\mbox{2.7e-17} & \mbox{-1.7e-16} & \mbox{9.7e-18} & \mbox{8.9e-23} &
\mbox{ 0.0647} & \mbox{ 0.0552} & 1.70
\\ 
 66 &  1156 &  1242 &  2246 &  4489 &  4488 &
\mbox{-1.0e-16} & \mbox{-1.0e-16} & \mbox{-1.0e-16} & \mbox{8.9e-23} &
\mbox{ 0.0616} & \mbox{ 0.0537} & 1.74
\\ 
 67 &  1190 &  1279 &  2313 &  4623 &  4623 &
\mbox{-1.1e-16} & \mbox{-2.2e-16} & \mbox{-1.7e-16} & \mbox{1.1e-22} &
\mbox{ 0.0609} & \mbox{ 0.0534} & 1.75
\\
\hline
\end{array}
\end{equation*}
\caption{Spherical $t$-designs on $\mathbb{S}^2$ with no symmetry restrictions, $N = \widehat{N}(2,t)$ and degrees $t = 61$ -- $120$}
\label{RSWtab:SF2}
\end{table}
\begin{table}[ht]
\begin{equation*}
\begin{array}{|c|ccccc|rrrc|ccc|}
\hline
 t & N^*(2,t) & N^+(2,t) & N & n & m & 
 \scriptstyle V_{t,N,\psi_1}(X_N) & \scriptstyle V_{t,N,\psi_2}(X_N) & \scriptstyle V_{t,N,\psi_3}(X_N) & \scriptstyle \bsr(X_N)^T\bsr(X_N) & 
 \delta(X_N) & h(X_N) & \rho(X_N)\\ 
\hline
121 &  3782 &  4089 &  7443 & 14883 & 14883 &
\mbox{-5.2e-16} & \mbox{-1.5e-16} & \mbox{-9.0e-16} & \mbox{3.3e-21} &
\mbox{ 0.0338} & \mbox{ 0.0295} & 1.75
\\ 
122 &  3844 &  4156 &  7566 & 15129 & 15128 &
\mbox{1.8e-16} & \mbox{1.7e-16} & \mbox{-8.9e-16} & \mbox{3.6e-21} &
\mbox{ 0.0356} & \mbox{ 0.0295} & 1.66
\\ 
123 &  3906 &  4224 &  7689 & 15375 & 15375 &
\mbox{1.8e-16} & \mbox{1.3e-16} & \mbox{-6.8e-16} & \mbox{4.0e-21} &
\mbox{ 0.0352} & \mbox{ 0.0292} & 1.66
\\ 
124 &  3969 &  4292 &  7814 & 15625 & 15624 &
\mbox{-9.1e-17} & \mbox{6.3e-17} & \mbox{-8.4e-16} & \mbox{4.1e-21} &
\mbox{ 0.0348} & \mbox{ 0.0298} & 1.71
\\ 
125 &  4032 &  4360 &  7939 & 15875 & 15875 &
\mbox{-8.7e-17} & \mbox{-1.9e-16} & \mbox{-1.0e-15} & \mbox{4.2e-21} &
\mbox{ 0.0337} & \mbox{ 0.0288} & 1.70
\\ 
126 &  4096 &  4430 &  8066 & 16129 & 16128 &
\mbox{-8.7e-17} & \mbox{-1.9e-16} & \mbox{-8.1e-16} & \mbox{4.4e-21} &
\mbox{ 0.0332} & \mbox{ 0.0283} & 1.71
\\ 
127 &  4160 &  4499 &  8193 & 16383 & 16383 &
\mbox{-7.6e-17} & \mbox{-6.8e-18} & \mbox{-8.1e-16} & \mbox{4.5e-21} &
\mbox{ 0.0321} & \mbox{ 0.0283} & 1.76
\\
\hline
\end{array}
\end{equation*}
\caption{Spherical $t$-designs on $\mathbb{S}^2$ with no symmetry restrictions, $N = \widehat{N}(2,t)$ and degrees $t = 121$ -- $180$}
\label{RSWtab:SF3}
\end{table}
\begin{table}[ht]
\begin{equation*}
\begin{array}{|c|ccccc|rrrc|ccc|}
\hline
 t & N^*(2,t) & N^+(2,t) & N & n & m & 
 \scriptstyle V_{t,N,\psi_1}(X_N) & \scriptstyle V_{t,N,\psi_2}(X_N) & \scriptstyle V_{t,N,\psi_3}(X_N) & \scriptstyle \bsr(X_N)^T\bsr(X_N) & 
 \delta(X_N) & h(X_N) & \rho(X_N)\\ 
\hline
  1 &     2 &     2 &     2 &     0 &     0 &
\mbox{0.0e+00} & \mbox{0.0e+00} & \mbox{0.0e+00} & \mbox{0.0e+00} &
\mbox{ 3.1416} & \mbox{ 1.5708} & 1.00
\\ 
  3 &     6 &     6 &     6 &     3 &     5 &
\mbox{3.5e-17} & \mbox{0.0e+00} & \mbox{0.0e+00} & \mbox{6.1e-31} &
\mbox{ 1.5708} & \mbox{ 0.9553} & 1.22
\\ 
  5 &    12 &    12 &    12 &     9 &    14 &
\mbox{6.9e-18} & \mbox{1.5e-17} & \mbox{-8.5e-17} & \mbox{5.2e-30} &
\mbox{ 1.1071} & \mbox{ 0.6524} & 1.18
\\ 
  7 &    20 &    20 &    32 &    29 &    27 &
\mbox{2.1e-17} & \mbox{4.0e-18} & \mbox{1.9e-16} & \mbox{3.4e-28} &
\mbox{ 0.5863} & \mbox{ 0.4480} & 1.53
\\ 
  9 &    30 &    31 &    48 &    45 &    44 &
\mbox{4.2e-18} & \mbox{-4.1e-18} & \mbox{-1.9e-16} & \mbox{1.2e-25} &
\mbox{ 0.4611} & \mbox{ 0.3860} & 1.67
\\ 
 11 &    42 &    43 &    70 &    67 &    65 &
\mbox{-2.4e-17} & \mbox{-1.9e-18} & \mbox{-1.4e-16} & \mbox{8.6e-27} &
\mbox{ 0.3794} & \mbox{ 0.3017} & 1.59
\\ 
 13 &    56 &    58 &    94 &    91 &    90 &
\mbox{-8.4e-18} & \mbox{4.8e-18} & \mbox{-1.3e-16} & \mbox{1.6e-25} &
\mbox{ 0.3146} & \mbox{ 0.2639} & 1.68
\\
\hline
\end{array}
\end{equation*}
\caption{Symmetric spherical $t$-designs on $\mathbb{S}^2$ with $N = \overline{N}(2,t)$ and odd degrees $t = 1$ -- $109$, except for $t = 1$ when $N = \overline{N}(2,t)-2$ and $t = 7, 11$ when $N = \overline{N}(2,t)+2$}
\label{RSWtab:SS1}
\end{table}
\begin{table}[ht]
\begin{equation*}
\begin{array}{|c|ccccc|rrrc|ccc|}
\hline
 t & N^*(2,t) & N^+(2,t) & N & n & m & 
 \scriptstyle V_{t,N,\psi_1}(X_N) & \scriptstyle V_{t,N,\psi_2}(X_N) & \scriptstyle V_{t,N,\psi_3}(X_N) & \scriptstyle \bsr(X_N)^T\bsr(X_N) & 
 \delta(X_N) & h(X_N) & \rho(X_N)\\ 
\hline
111 &  3192 &  3449 &  6218 &  6215 &  6215 &
\mbox{1.3e-16} & \mbox{-1.9e-16} & \mbox{-6.7e-16} & \mbox{6.6e-21} &
\mbox{ 0.0376} & \mbox{ 0.0324} & 1.72
\\ 
113 &  3306 &  3573 &  6444 &  6441 &  6440 &
\mbox{1.1e-16} & \mbox{-1.5e-16} & \mbox{-7.2e-16} & \mbox{7.2e-21} &
\mbox{ 0.0377} & \mbox{ 0.0318} & 1.69
\\ 
115 &  3422 &  3698 &  6672 &  6669 &  6669 &
\mbox{2.2e-16} & \mbox{2.0e-16} & \mbox{-6.7e-16} & \mbox{8.8e-21} &
\mbox{ 0.0375} & \mbox{ 0.0316} & 1.69
\\ 
117 &  3540 &  3826 &  6906 &  6903 &  6902 &
\mbox{-3.7e-16} & \mbox{1.6e-17} & \mbox{-5.6e-16} & \mbox{8.6e-21} &
\mbox{ 0.0358} & \mbox{ 0.0309} & 1.73
\\ 
119 &  3660 &  3957 &  7142 &  7139 &  7139 &
\mbox{-4.0e-16} & \mbox{2.0e-16} & \mbox{-9.7e-16} & \mbox{9.5e-21} &
\mbox{ 0.0356} & \mbox{ 0.0303} & 1.71
\\ 
121 &  3782 &  4089 &  7384 &  7381 &  7380 &
\mbox{-5.0e-16} & \mbox{-9.3e-17} & \mbox{-8.7e-16} & \mbox{1.1e-20} &
\mbox{ 0.0351} & \mbox{ 0.0297} & 1.69
\\ 
123 &  3906 &  4224 &  7628 &  7625 &  7625 &
\mbox{1.7e-16} & \mbox{1.3e-16} & \mbox{-7.7e-16} & \mbox{1.2e-20} &
\mbox{ 0.0349} & \mbox{ 0.0292} & 1.67
\\
\hline
\end{array}
\end{equation*}
\caption{Symmetric spherical $t$-designs on $\mathbb{S}^2$ with $N = \overline{N}(2,t)$ and odd degrees $t = 111$ -- $219$}
\label{RSWtab:SS2}
\end{table}
\begin{table}[ht]
\begin{equation*}
\begin{array}{|c|ccccc|rrrc|ccc|}
\hline
 t & N^*(2,t) & N^+(2,t) & N & n & m & 
 \scriptstyle V_{t,N,\psi_1}(X_N) & \scriptstyle V_{t,N,\psi_2}(X_N) & \scriptstyle V_{t,N,\psi_3}(X_N) & \scriptstyle \bsr(X_N)^T\bsr(X_N) & 
 \delta(X_N) & h(X_N) & \rho(X_N)\\ 
\hline
221 & 12432 & 13488 & 24534 & 24531 & 24530 &
\mbox{-1.4e-15} & \mbox{2.6e-16} & \mbox{-1.4e-15} & \mbox{4.0e-19} &
\mbox{ 0.0194} & \mbox{ 0.0163} & 1.68
\\ 
223 & 12656 & 13732 & 24978 & 24975 & 24975 &
\mbox{-8.1e-17} & \mbox{-1.6e-16} & \mbox{-1.3e-15} & \mbox{4.4e-19} &
\mbox{ 0.0192} & \mbox{ 0.0164} & 1.71
\\ 
225 & 12882 & 13978 & 25428 & 25425 & 25424 &
\mbox{-4.9e-16} & \mbox{-5.3e-16} & \mbox{-1.3e-15} & \mbox{4.5e-19} &
\mbox{ 0.0188} & \mbox{ 0.0160} & 1.70
\\ 
227 & 13110 & 14226 & 25880 & 25877 & 25877 &
\mbox{1.5e-16} & \mbox{6.4e-16} & \mbox{-1.4e-15} & \mbox{4.7e-19} &
\mbox{ 0.0185} & \mbox{ 0.0158} & 1.72
\\ 
229 & 13340 & 14476 & 26338 & 26335 & 26334 &
\mbox{-1.0e-15} & \mbox{1.8e-16} & \mbox{-1.5e-15} & \mbox{5.1e-19} &
\mbox{ 0.0183} & \mbox{ 0.0156} & 1.71
\\ 
231 & 13572 & 14728 & 26798 & 26795 & 26795 &
\mbox{3.6e-17} & \mbox{-5.2e-16} & \mbox{-1.4e-15} & \mbox{5.3e-19} &
\mbox{ 0.0181} & \mbox{ 0.0157} & 1.74
\\ 
233 & 13806 & 14982 & 27264 & 27261 & 27260 &
\mbox{-6.2e-16} & \mbox{-2.4e-16} & \mbox{-1.5e-15} & \mbox{5.7e-19} &
\mbox{ 0.0185} & \mbox{ 0.0155} & 1.68
\\
\hline
\end{array}
\end{equation*}
\caption{Symmetric spherical $t$-designs on $\mathbb{S}^2$ with $N = \overline{N}(2,t)$ and odd degrees $t = 221$ -- $325$}
\label{RSWtab:SS3}
\end{table}
\begin{table}[ht]
\begin{equation*}
\begin{array}{|c|ccccc|rrr|rrc|}
\hline
 t & N^*(3,t) & N^+(3,t) & \widehat{N}(3,t) & N & sym & 
 V_{t,N,\psi_1}(X_N) & V_{t,N,\psi_2}(X_N) & V_{t,N,\psi_3}(X_N) & 
 \delta(X_N) & h(X_N) & \rho(X_N)\\ 
\hline
  1 &     2 &     2 &     4 &     2 & 1 &
\mbox{5.6e-17} & \mbox{0.0e+00} & \mbox{0.0e+00} &
\mbox{ 3.1416} & \mbox{ 1.5708} & 1.00
\\ 
  2 &     5 &     5 &     7 &     5 & 0 &
\mbox{8.4e-17} & \mbox{4.4e-17} & \mbox{1.7e-16} &
\mbox{ 1.8235} & \mbox{ 1.3181} & 1.45
\\ 
  3 &     8 &     8 &    12 &     8 & 1 &
\mbox{1.5e-17} & \mbox{-1.1e-16} & \mbox{-2.6e-17} &
\mbox{ 1.5708} & \mbox{ 1.0472} & 1.33
\\ 
  3 &     8 &     8 &    12 &    16 & 1 &
\mbox{-1.4e-17} & \mbox{-9.6e-17} & \mbox{-1.5e-16} &
\mbox{ 1.0472} & \mbox{ 1.0472} & 2.00
\\ 
  5 &    20 &    19 &    32 &    24 & 1 &
\mbox{-4.7e-17} & \mbox{1.3e-16} & \mbox{-6.8e-17} &
\mbox{ 1.0472} & \mbox{ 0.7854} & 1.50
\\ 
  7 &    40 &    40 &    70 &    48 & 1 &
\mbox{2.8e-17} & \mbox{-1.2e-16} & \mbox{1.8e-16} &
\mbox{ 0.7854} & \mbox{ 0.6086} & 1.55
\\ 
 11 &   112 &   117 &   219 &   120 & 1 &
\mbox{-1.7e-17} & \mbox{-5.6e-17} & \mbox{-5.7e-16} &
\mbox{ 0.6283} & \mbox{ 0.3881} & 1.24
\\ 
 11 &   112 &   117 &   219 &   600 & 1 &
\mbox{2.9e-17} & \mbox{-3.7e-17} & \mbox{-8.1e-17} &
\mbox{ 0.2709} & \mbox{ 0.3881} & 2.87
\\ 
\hline
\end{array}
\end{equation*}
\caption{Regular spherical $t$-designs on $\mathbb{S}^3$ for degrees $t$ = 1, 2, 3, 5 ,7 ,11}
\label{RSWtab:S3SDR}
\end{table}
\begin{table}[ht]
\begin{equation*}
\begin{array}{|c|ccccc|rrr|rrc|}
\hline
 t & N^*(3,t) & N^+(3,t) & N & n & m & 
 V_{t,N,\psi_1}(X_N) & V_{t,N,\psi_2}(X_N) & V_{t,N,\psi_3}(X_N) & 
 \delta(X_N) & h(X_N) & \rho(X_N)\\ 
\hline
  1 &     2 &     2 &     4 &     6 &   4 & \mbox{6.2e-17} & \mbox{0.0e+00} & \mbox{0.0e+00} &
\mbox{ 1.5708} & \mbox{ 1.5708} & 2.00
\\ 
  2 &     5 &     5 &     7 &    15 &  13 & \mbox{6.3e-17} & \mbox{5.6e-17} & \mbox{1.0e-16} &
\mbox{ 1.3585} & \mbox{ 1.1683} & 1.72
\\ 
  3 &     8 &     8 &    12 &    30 &  29 & \mbox{8.1e-18} & \mbox{-9.1e-17} & \mbox{-7.9e-17} &
\mbox{ 1.2311} & \mbox{ 1.0016} & 1.63
\\ 
  4 &    14 &    13 &    20 &    54 &  54 & \mbox{1.4e-19} & \mbox{5.6e-17} & \mbox{-2.4e-17} &
\mbox{ 0.9414} & \mbox{ 0.8338} & 1.77
\\ 
  5 &    20 &    19 &    32 &    90 &  90 & \mbox{-3.0e-17} & \mbox{1.3e-16} & \mbox{2.6e-17} &
\mbox{ 0.7816} & \mbox{ 0.7041} & 1.80
\\ 
  6 &    30 &    28 &    49 &   141 & 139 & \mbox{1.1e-17} & \mbox{1.0e-17} & \mbox{-9.2e-18} &
\mbox{ 0.6883} & \mbox{ 0.6086} & 1.77
\\ 
  7 &    40 &    40 &    70 &   204 & 203 & \mbox{2.1e-17} & \mbox{-1.2e-16} & \mbox{-1.8e-16} &
\mbox{ 0.5765} & \mbox{ 0.5454} & 1.89
\\
\hline
\end{array}
\end{equation*}
\caption{Computed spherical $t$-designs on $\mathbb{S}^3$ for degrees $t$ = 1,...,20, with $N = \widehat{N}(3,t)$}
\label{RSWtab:S3SDC}
\end{table}
\begin{table}[ht]
\begin{equation*}
\begin{array}{|c|ccccc|rrr|rrc|}
\hline
 t & N^*(3,t) & N^+(3,t) & N & n & m & 
 V_{t,N,\psi_1}(X_N) & V_{t,N,\psi_2}(X_N) & V_{t,N,\psi_3}(X_N) & 
 \delta(X_N) & h(X_N) & \rho(X_N)\\ 
\hline
  1 &     2 &     2 &     4 &     6 &   4 & \mbox{6.2e-17} & \mbox{0.0e+00} & \mbox{0.0e+00} &
\mbox{ 1.5708} & \mbox{ 1.5708} & 2.00
\\ 
  3 &     8 &     8 &    10 &    30 &  29 & \mbox{-2.4e-17} & \mbox{-1.1e-16} & \mbox{-2.6e-16} &
\mbox{ 1.3181} & \mbox{ 0.9776} & 1.48
\\ 
  5 &    20 &    19 &    28 &    90 &  90 & \mbox{-2.0e-17} & \mbox{1.3e-16} & \mbox{9.2e-17} &
\mbox{ 0.8334} & \mbox{ 0.7303} & 1.75
\\ 
  7 &    40 &    40 &    60 &   204 & 203 & \mbox{5.7e-18} & \mbox{-1.3e-16} & \mbox{-1.1e-16} &
\mbox{ 0.6324} & \mbox{ 0.5656} & 1.79
\\ 
  9 &    70 &    71 &   114 &   384 & 384 & \mbox{-2.7e-17} & \mbox{6.4e-17} & \mbox{-1.4e-16} &
\mbox{ 0.4863} & \mbox{ 0.4548} & 1.87
\\ 
 11 &   112 &   117 &   194 &   651 & 649 & \mbox{1.7e-17} & \mbox{-4.7e-17} & \mbox{-2.2e-16} &
\mbox{ 0.4126} & \mbox{ 0.3860} & 1.87
\\ 
 13 &   168 &   178 &   308 &  1014 & 1014 & \mbox{3.3e-17} & \mbox{5.5e-17} & \mbox{-1.1e-16} &
\mbox{ 0.3454} & \mbox{ 0.3220} & 1.87
\\
\hline
\end{array}
\end{equation*}
\caption{Computed symmetric spherical $t$-designs on $\mathbb{S}^3$ for degrees $t$ = 1,3,...,31, with $N = \overline{N}(3,t)$}
\label{RSWtab:S3SDS}
\end{table}

\begin{acknowledgement}
This research includes extensive computations using the Linux computational cluster Katana
supported by the Faculty of Science, UNSW Sydney.
\end{acknowledgement}

%
%


\end{document}